\documentclass[12pt,a4paper]{article}
\usepackage{amsfonts,amssymb,amsmath}
\usepackage[english]{babel}
\usepackage{amssymb}
\usepackage{amsthm}
\usepackage{pdfsync}
 \usepackage{color,makeidx}
\usepackage{hyperref}
\usepackage{mathrsfs}
 \usepackage{enumitem}

-\marginparwidth \oddsidemargin -4mm \evensidemargin \oddsidemargin

\newcommand{\eps}{\varepsilon}
\newcommand{\commentout}[1]{}
\newcommand{\E}{{\mathbb E}}
\newcommand{\bE}{\mathbb{E}}

\newcommand{\bbE}{\mathbb{E}}
\newcommand{\bbP}{\mathbb P}

\newcommand{\bbZ}{\mathbb Z}

\newcommand{\ii}{{\rm i}}

\numberwithin{equation}{section}
\theoremstyle{plain}
\newtheorem{thm}{Theorem}[section]
\newtheorem{lemma}[thm]{Lemma}
\newtheorem{proposition}[thm]{Proposition}
\newtheorem{corollary}[thm]{Corollary}

\theoremstyle{definition}

\newtheorem{remark}[thm]{Remark}

\newcommand{\bes}{\begin{displaymath}}
\newcommand{\ees}{\end{displaymath}}
\newcommand{\be}{\begin{equation}}
\newcommand{\ee}{\end{equation}}
\newcommand{\ba}{\begin{eqnarray}}
\newcommand{\ea}{\end{eqnarray}}
\newcommand{\bas}{\begin{eqnarray*}}
\newcommand{\eas}{\end{eqnarray*}}

\newcommand{\C}{{\@Bbb C}}

\renewcommand{\P}{{\@Bbb P}}

\newcommand{\Q}{{\@Bbb Q}}
\newcommand{\bQ}{{\@Bbb Q}}
\newcommand{\Z}{{\mathbb{Z}}}
\newcommand{\al}{\alpha}
\newcommand{\si}{\sigma}

\newcommand{\Om}{\Omega}

\newcommand{\bze}{ 0}

\newcommand{\cA}{\@s A}
\newcommand{\cB}{\@s B}
\newcommand{\cC}{\@s C}
\newcommand{\cD}{\@s D}
\newcommand{\cE}{\@s E}
\newcommand{\cF}{\@s F}
\newcommand{\cG}{\@s G}
\newcommand{\cH}{\@s H}
\newcommand{\cI}{\@s I}
\newcommand{\cJ}{\@s J}
\newcommand{\cK}{\@s K}
\newcommand{\cL}{\@s L}
\newcommand{\cN}{\@s N}
\newcommand{\cM}{\@s M}
\newcommand{\cO}{\@s O}
\newcommand{\cP}{\@s P}
\newcommand{\cR}{\@s R}
\newcommand{\cS}{\@s S}
\newcommand{\cT}{\@s T}
\newcommand{\cV}{\@s V}
\newcommand{\cW}{\@s W}
\newcommand{\cX}{\@s X}
\newcommand{\cY}{\@s Y}
\newcommand{\cZ}{\@s Z}

\newcommand{\bma}{\@bm a}
\newcommand{\bmb}{\@bm b}
\newcommand{\bmc}{\@bm c}
\newcommand{\bmd}{\@bm d}

\newcommand{\bme}{\@bm e}
\newcommand{\bmf}{\@bm f}
\newcommand{\bmg}{\@bm g}
\newcommand{\bmh}{\@bm h}
\newcommand{\bmi}{\@bm i}
\newcommand{\bmj}{\@bm j}
\newcommand{\bmk}{\@bm k}
\newcommand{\bml}{\@bm l}
\newcommand{\bmm}{\@bm m}
\newcommand{\bmn}{\@bm n}
\newcommand{\bmo}{\@bm o}
\newcommand{\bmp}{\@bm p}
\newcommand{\bmq}{\@bm q}
\newcommand{\bmr}{\@bm r}
\newcommand{\bms}{\@bm s}
\newcommand{\bmt}{\@bm t}
\newcommand{\bmu}{\@bm u}
\newcommand{\bmw}{\@bm w}
\newcommand{\bmv}{\@bm v}
\newcommand{\bmx}{\@bm x}
\newcommand{\bx}{\@bm x}
\newcommand{\bmy}{\@bm y}
\newcommand{\bmz}{\@bm z}
\newcommand{\bbN}{\mathbb N}
\newcommand{\by}{\@bm y}
\newcommand{\bmzero}{\@bm 0}

\newcommand{\ga}{\gamma}

\newcommand{\gA}{\@g A}
\newcommand{\gD}{\@g D}
\newcommand{\gJ}{\@g J}
\newcommand{\gF}{\@g F}
\newcommand{\gM}{\@g M}
\newcommand{\gR}{\@g R}
\newcommand{\bbR}{\mathbb R}

\newcommand{\varcal}[1]{\mathscr{#1}}
\newcommand{\R}{\mathbb R}

 \newcommand{\N}{\mathbb N}
\newcommand{\Rd}{{\mathbb R^d}}
\newcommand{\Rdp}{{\mathbb R^{d+1}}}
\newcommand{\F}{\varcal{F}}

\renewcommand{\P}{\mathbb P}

\newcommand{\SprimeRd}{%
 {\ifmmode {\varcal S}'(\Rd)%
 \else ${\varcal S}'(\Rd)$\fi}}
\newcommand{\SRd}{{%
 \ifmmode \varcal S(\Rd)%
 \else $\varcal S(\Rd)$\fi}}
\newcommand{\Sprime}{{%
 \ifmmode {\varcal S}'%
 \else ${\varcal S}'$\fi}}
\newcommand{\Ss}{{%
 \ifmmode \varcal S%
 \else $\varcal S$\fi}}
\newcommand{\G}{\varcal{G}}

\newcommand{\D}{\varcal{D}}
\def\ind{1\mkern-7mu1}
\long\def\arr>#1>#2>{\xrightarrow[#2]{#1}}
\long\def\larr<#1<#2<{\xleftarrow[#2]{#1}}

\newcommand{\abs}[1]{\left\vert#1\right\vert} 

\newcommand{\bfnorm}[1]{\vert \kern -1.9pt \vert #1 \vert \kern -1.9pt \vert}

\newcommand{\W}{U} 

\title{Non-Gaussian limit of a tracer motion in an  incompressible flow}

\begin{document}

\author{Tomasz Komorowski \thanks{IMPAN,
ul. \'{S}niadeckich 8,   00-956 Warsaw, Poland, e-mail:
komorow@hektor.umcs.lublin.pl}, Anna Talarczyk \thanks{IMPAN,
ul. \'{S}niadeckich 8,   00-956 Warsaw, Poland and Institute of Mathematics, University of Warsaw, ul. Banacha 2, 02-097 Warsaw, Poland  e-mail:
annatal@mimuw.edu.pl}}

\maketitle
\begin{abstract}
We consider a massless tracer particle moving in a random, stationary, isotropic and divergence free velocity field. We identify a class of  fields, for which the limit of the laws of appropriately scaled tracer trajectory processes is non-Gaussian but a Rosenblatt type  process.
\end{abstract}

\textbf{}


\medskip


\section{Introduction}

We consider the motion of a massless tracer particle in a random field described by an ordinary differential equation
\begin{equation}
\label{011706}
\frac{dx(t)}{dt}=V(t,x(t)),\quad x(0)=x_0.
\end{equation}
Here $V=(V_1,\ldots,V_d):\bbR^{1+d}\times\Om\to\bbR^d$ is a $d$-dimensional random vector field defined over a probability space $(\Om,{\cal V},\bbP)$.
It is one of the most popular models of transport in turbulence
considered in statistical hydrodynamics that can be traced to the
works of  G. Taylor in the 1920's, see \cite{taylor} and R. Kraichnan
in the 1970's, see \cite{kraichnan}. There exists an extensive 
literature concerning the model, see e.g. \cite{frisch,gawedzki,majda-kramer,yaglom} and the references therein. Due to this motivation it is often assumed that the field $V(t,x)$,  thought of as the Eulerian velocity of the fluid, is incompressible, i.e. 
$$
\nabla_x\cdot V(t,x):=\sum_{j=1}^d\partial_{x_j}V_j(t,x)\equiv 0.
$$
Another assumption frequently made about the velocity flow is its
stationarity and statistical isotropy. The above means that the law of
the random field $\left(V(t,x)\right)_{(t,x)\in\bbR^{1+d}}$ does not depend on both
temporal and spatial translations and  is also  invariant under the action
of the group
of  rotations in space. As a consequence of the incompressibiliy and
stationarity of the
flow one concludes that the velocity observed along the
trajectory $\left(V(t,x(t)\right)_{t\in\bbR}$  - the so called {\em
  Lagrangian process} -  is stationary, see Theorem 2, p. 500 of \cite{port-stone}.

When the field is of zero mean and its fluctuations are small, i.e. \eqref{011706} can be written in the form
\begin{equation}
\label{011706a}
\frac{dx(t)}{dt}=\eps V(t,x(t)),\quad x(0)=x_0,
\end{equation}
where $\eps>0$ is a certain parameter (that turns out to be small),
then the scaled trajectory of the tracer $x(t\eps^{-2})$ satisfies the
central limit theorem (CLT), as $\eps\to0+$, provided that the field
is sufficiently strongly mixing in the temporal variable,  see
\cite{khasminskii} and also \cite{kp}. More precisely, one can show
the functional CLT, {i.e.  the continuous trajectory} processes  $\left(x(t\eps^{-2})\right)_{t\ge0}$ converge in law, as $\eps\to0+$, to a zero mean Brownian motion whose covariance is given by 
 the Taylor-Kubo formula
\be
\label{TK}
D_{j,j'}:=\int\limits^{+\infty}_0\left\{\bE[V_j(t,\bze)V_{j'}(0,\bze)]+
\bE[V_{j'}(t,\bze)V_{j}(0,\bze)]\right\}dt,\quad j,j'=1,\ldots,d.
\ee

Consider now a  stationary and isotropic, divergence free  field 
whose  energy spectrum satisfies the power law. Its covariance matrix can be then  written  as
\begin{equation}
\label{010903}
R_{j,j'}(t,x)=\bbE\left[V_j(t+t',x+x')V_{j'}(t',x')\right]=
\int_{\bbR^{d+1}}\exp\left\{\ii k\cdot x\right\}  \Gamma_{j,j'}(k) \frac{\hat{\cal E}\left(t,|k|\right)dk }{|k|^{d-1}},\quad
\end{equation}
for any  $j,j'=1,\ldots,d$, $t,t'\in\bbR$, $x,x'\in\bbR^d$.   Here $\ii=\sqrt{-1}$ and
\begin{equation}
\label{etk}
\hat{\cal E}\left(t,\xi\right)=e^{-r(\xi)|t|}a(\xi){\cal E}_0(\xi),\quad (t,\xi)\in\bbR\times\bbR_+.
\end{equation}
 We assume that the cut-off
function $a(\cdot)$ is continuous, non-negative, compactly supported
with $a(0)>0$.   Parameter $r(\cdot)$, determining the mixing
rate in time, and the power-energy density spectrum ${\cal E}_0(\cdot)$  are both non-negative and satisfy $r(\xi)\sim
\xi^{2\beta}$ and 
  ${\cal E}_0(\xi)\sim \xi^{1-2\al}$, as $\xi\ll1$ for
some parametrs $\al,\beta$. The factor $\Gamma_{j,j'}(k):=
\delta_{j,j'}-k_jk_{j'}/|k|^2$ ensures that the realizations of the
field are incompressible.
We can see from  formula \eqref{010903} that the spectrum is integrable in
$k$ and the mixing rate decays on large spatial scales, provided that   $\al<1$ and $\beta\ge0$.

Random fields whose covariance  is given by \eqref{010903}  play an important role in statistical
hydrodynamics. According to the theory of Kolmogorov and Obukhov, see
\cite{kolm,obuk}, the Eulerian velocity of a
fully turbulent flow of a fluid can be described by a locally
isotropic, time-space homogenous random field whose energy spectrum is
of the aforementined form. 

For  a random field whose covariance is given by \eqref{010903}         
the coefficients $D_{j,j'}$ appearing in Taylor-Kubo formula \eqref{TK} are finite if and only if
\begin{equation}
\label{C1}
\int_0^{+\infty}\frac{a(\xi)d\xi}{\xi^{2\al+2\beta-1}}<+\infty,
\end{equation}
which leads to the condition $\al+\beta<1$.
For  a tracer motion in a Gaussian field whose  covariance satisfies \eqref{C1} the functional CLT has been proved in \cite{fannjiang-kom}. 
One can inquire what would happen when the condition \eqref{C1} is
violated. The conjecture that the trajectory evolves on a shorter,
superdiffusive scale turns out to be true, see
\cite{fannjiang-kom2}. Namely, it has been shown that if
$\al+\beta>1$ (and $\al<1$), then the scaled processes
$\left(x(t\eps^{-2\delta})\right)$, with
$\delta:=\beta/(\al+2\beta-1)$, converge in law, as $\eps\to0+$, to a zero mean fractional Brownian motion with Hurst exponent $H=1/(2\delta)$.

A natural question arises whether a fractional Brownian motion is a universal limit for an anomalous diffusive scaling of the tracer particle. More precisely,  take  a family of non-Gaussian fields with the covariance
of the form \eqref{010903}
 and consider the situation when $\al+\beta>1$, so the diffusivity given  the Kubo formula  is no longer finite.
Does the limit of an appropriately scaled trajectory $\left(x(t\eps^{-2\delta})\right)$ tend to a Gaussian, as $\eps\to0+$? 
The experimental data in general seem to confirm Gaussianity of the tracer distribution, see e.g. \cite{hu}, but the non-Gaussianity is also observed in some explicitly solvable models, see e.g. Sec 5.2.2 of  \cite{majda-kramer}. We are not aware of any rigorous result for isotropic flows that  leads to a non-Gaussian limit.

In the present paper we consider fields  whose components belong to the space ${\cal H}_2$ in the Fock's decomposition over {some} Gaussian
Hilbert space. We show that then the limit of the tracer particle is
a stochastic process that is related to the  Rosenblatt  one, see Theorem \ref{thm:main},  Remark \ref{rem:2} and Proposition \ref{prop:Rosenblatt}
below. 
The Rosenblatt process has been introduced by Taqqu in \cite{Taqqu1975}. It is self-similar, with stationary increments and ``lives'' in the 
 second Wiener chaos (see \eqref{e:Rosenblatt} for one of its possible
 representations). {Unlike in  the  Gaussian case, in the second
 chaos there exist many self-similar processes with stationary
 increments (see \cite{MT}), but the Rosenblatt one} is the simplest of
 them and is considered to be the second chaos counterpart of the
 fractional Brownian motion. This process, along with other self
 similar processes with stationary increments living in the second and
 higher chaoses, {have attracted 
   quite a lot of interest lately}. 
{Some recent literature concerning the subject includes (but
certainly not exhausts)} \cite{T4,tudor, MT,MT1,Arras,BGT_Rosenblatt,BT1,BT2}.
 Our limit is not strictly speaking a Rosenblatt process, but it  is
 closely related to a process of that type, {see
   formula \eqref{e:10.17} and Proposition \ref{prop:Rosenblatt} below.}
 As a by product  we obtain a representation of the Rosenblatt process that seems to be new, according to our knowledge, see Proposition \ref{prop:Rosenblatt}.


{Concerning the proof of our main result, see  Theorem \ref{thm:main}
below, we compute the limits of the moments of the
trajectory process, see Theorem \ref{thm090801} below. In fact, for convenience sake we consider the
moments of  its statistically equivalent version introduced in
\eqref{e:10.24} below. 
These limits   coincide with the corresponding
moments of the limiting process in question. This, in turn implies the
convergence of 
finite dimensional laws, as the respective moment problem for the  limiting
law is well posed.
The proof of Theorem \ref{thm090801}  relies on showing that the limits of the moments of the trajectory process 
coincide with the limits of the respective moments 
of the approximate process
 $\eps\int_0^{t/\eps^{2\delta}}V(s,x_0)ds$, obtained by ''freezing'' the right
 hand side at the initial position of the tracer. The convergence of
 the latter is quite straightforward due to Proposition \ref{lem:10.6}. }
{To argue that the ''true'' trajectory  can be replaced by its
approximate counterpart  we use the Taylor expansion
\eqref{expansion1}. An important role in  the analysis of the arising
terms  is played by the representation of 
the  products of multiple
stochastic integrals in terms of integrals of higher order with the
help of appropriate Feynman diagrams (Proposition \ref{lem:A1}) and
the formula for their conditional moments (Proposition \ref{lem:A3}).}

The organization of the paper is as follows: in Section \ref{sec2} we
introduce basic notions and formulate our main result, see Theorem
\ref{thm:main}. In Section \ref{sec3a} we recall some basic facts
about multiple stochastic integration in the context of Gaussian
stochastic measures. Section \ref{sec3} is devoted to the proof of the
main result. Finally in Section \ref{sec5} we prove 
Proposition \ref{prop:Rosenblatt} concerning a representation of the
Rosenblatt process.

\subsection*{Acknowledgements} The research of T. K. has been
supported by  National Science Center, Poland, grant \linebreak 2016/23/B/ST1/00492.
The research of A. T. has been supported by  National Science Center,
Poland, grant  DEC-2012/07/B/ST1/03417.

\section{Preliminaries and formulation of the main result}

\label{sec2}

\subsection{Some basic notation}

We shall use the following notation: for any two functions
$f,g:A\to\bbR$ and some cone $A\subset\bbR^d$ we write
$$
f(x)\preceq g(x),\quad x\in A
$$
iff there exist $C,c>0$ such that
$$
f(x)\le Cg(c x),\quad x\in A.
$$
We write 
$$
f(x)\approx g(x),\quad x\in A
$$
iff
$$
f(x)\preceq g(x)\quad \mbox{and}\quad g(x)\preceq f(x),\quad x\in A.
$$
We shall also denote by $\bbZ_+\,(\bbN)$,  the set of all non-negative
(positive) integers, $\bbZ_d:=\{1,\ldots,d\}$ for some positive
integer $d$
and $\R_+$, $\bar \R_+$ stand for $(0,+\infty)$ and $[0,\infty)$, respectively.

Throughout the paper $C$, $C_i$ always denote positive constants which may be different from line to line. Possible dependencies on parameters is indicated in parenthesis.

With boldface we denote vectors, e.g. ${\bf t}=(t_1,\ldots, t_n)$, to stress the dimension we also often write ${\bf t}_{1,n}$, we also write ${\bf t}_{m,n}=(t_m,\ldots, t_n)$ for $m\le n$.
The differentials $ds_m\ldots ds_n$ for $m\le n$ are abbreviated as $d{\bf s}_{m,n}$.

$\Delta_n$ stands for the $n$ dimensional simplex 
$$
\Delta_n=[(s_1,\ldots, s_n)\in \R_+^n:s_1\ge s_2\ge\ldots \ge s_n\ge
0].
$$
{We shall also frequently use sets of the form
\begin{align}
\label{010810}
&\Box(\textbf{t}_{1,m}):= [0,t_1]\times\ldots\times
  [0,t_m],\quad\Delta(\textbf{t}_{1,m}):=\Delta_m\cap \Box(\textbf{t}_{1,m}),\nonumber\\
&\Delta_N({\bf t}_1,_m):=\Delta_{m+N}\cap\left(\Box(\textbf{t}_{1,m})\times\R^{N}\right),\\
&\Delta_N(s',s'') :=[(s_1,\ldots,s_N):\,s''\ge s_1\ge\ldots\ge
s_N\ge s'] 
\mbox{ for }0\le s'\le s'',\nonumber\\
&\Delta_{m}({\bf t_{1,m}}, s)=[(s_1,\ldots, s_m)\in
  \Delta(\textbf{t}_{1,m}): s_m\ge s]\mbox{ for }0\le s.\nonumber
\end{align}}

\subsection{Space-time white noise field}

\label{sec2.2}
Denote by $L^2(\bbR^{1+d})$ the complex Hilbert space of all
$\phi:\bbR^{1+d}\to\mathbb C$, for which
$$
\|\phi\|_{L^2(\bbR^{1+d})}^2=\int_{\bbR^{1+d}}|\phi(t,k)|^2dtdk<+\infty.
$$
The scalar product on $L^2(\bbR^{1+d})$ will be denoted by 
$\langle \cdot , \cdot\rangle_{L^2(\bbR^{1+d})}$.

We let
$L^2_{(s)}(\bbR^{1+d})$ be the real Hilbert space consisting of 
$\phi\in L^2(\bbR^{1+d})$   that satisfy $\phi^*(t,k)=\phi(t,-k)$. 
Here $z^*$ denotes the complex conjugate of $z\in\mathbb C$. 
The scalar product on $L^2_{(s)}(\bbR^{1+d})$ is given by a formula
$$
\langle \psi , \phi\rangle_{L^2(\bbR^{1+d})}=\int_{\bbR^{1+d}}\psi(t,k)\phi(t,-k)dtdk,\quad
\psi ,\phi\in L^2_{(s)}(\bbR^{1+d}).
$$

Let  $(e_m)_{m\ge1}$ be an orthonormal
base in $L^2_{(s)}(\bbR^{1+d})$.
A space-time, {$d$-dimensional vector valued}, white noise  $W(dt,dk)=(W_1(dt,dk),\ldots, W_d(dt,dk))$, $(t,k)\in \bbR^{1+d}$
can be defined as an $\bbR^d$-valued stochastic measure,  over a probability
space $(\Om,{\cal V},\bbP)$, in the following way: 
for any $\phi\in L^2_{(s)}(\bbR^{1+d})$ we let
\begin{equation}
\label{e:Wj}
\langle \phi,W_j\rangle:=\sum_{m=1}^{+\infty}\xi_{j,m}\langle \phi,e_m\rangle_{L^2(\bbR^{1+d})},\quad j=1,\ldots,d
\end{equation}
where $\left(\xi_{j,m}\right)$ are i.i.d. one dimensional, real
valued, standard normal random variables. The series on the right hand side of 
\eqref{e:Wj}  converges both 
a.s. and in the $L^2$ sense.
Observe that, by the  Parseval identity
\begin{equation}
\label{010403}
\bbE\left[\langle \psi ,W_j\rangle\langle \phi,W_{j'}\rangle \right]=\delta_{j,j'}\langle \psi , \phi\rangle_{L^2(\bbR^{1+d})},\quad j,j'=1,\ldots,d,\quad
\psi ,\phi\in L^2_{(s)}(\bbR^{1+d}),
\end{equation}
where $\delta_{j,j'}$ denotes the Kronecker delta function.

\subsection{Velocity field}\label{sec:velocity_field}
{In the present section we introduce the random fields considered
throughout the paper.}  The assumptions {made} in this section are valid throughout the paper.

 Suppose that   $a:\bar\bbR_+\to\bar \bbR_+$ is
 compactly supported, continuous function
satisfying
\begin{equation}
\label{Erb}
a_0:=a(0)>0.
\end{equation}
We shall also assume that $r:\bbR_+\to\bbR_+$ is a continuous function that satisfies
\begin{equation}
\label{022806a}
r(\xi)\approx \xi^{2\beta},\quad \xi\in\bbR_+\quad \text{and}\quad \lim_{\xi\to0}\frac{r(\xi)}{\xi^{2\beta}}=r_0
\end{equation}
for some $\beta>0$, $r_0>0$. We also fix a parameter $\alpha\in \R$, $\alpha<1$.

 Consider a random  vector field $V(t,x)=(V_1(t,x),\ldots,V_d(t,x))$ that is formed over the space of the second degree Hermite polynomials corresponding to the Gaussian Hilbert space for the noise $W(t,k)$, see e.g. Chapter 2 of \cite{janson}.  More precisely, for any $j=1,\ldots,d$ we let 
 \begin{align}
V_{j}(t,x):=\int_{(-\infty,t]^2}\int_{\bbR^{2d}}\exp\left\{\ii(k+k')\cdot
  x\right\}E(t-s,t-s',k,k') Y_{j}(ds,dk,ds',dk'),
\label{e:10.1}
\end{align}
where
\begin{equation}
\label{E}
E(s,s',k,k'):=\frac{\sqrt{r(\abs k ){a(\abs k)}}}{|k|^{(d+\alpha -1)/2}}\frac{\sqrt {r(\abs {k'})a(\abs {k'})}}{|k'|^{(d+\alpha -1)/2}}
\exp\left\{-\frac12\left[r(|k|)|s|+r(|k'|)|s'|\right]\right\}
\end{equation}
 \begin{align}
 \label{Y}
Y_{j}(ds,dk,ds',dk'):=&\sum_{j'=1}^d\Gamma_{j,j'}(k+k') W_{j'}(ds,dk)W_{j'}(ds',dk'),\\
 \Gamma_{j,j'}(k):=&
\delta_{j,j'}-k_jk_{j'}/|k|^2.
\label{e:Gamma}
\end{align}
Matrix $[\Gamma_{j,j'}(k) ]$ corresponds to the orthogonal
projection of $\bbR^d$ onto the hyperplane orthogonal to $k$, therefore  
\begin{equation}
 \label{e:Gamma2}
\sum_{\ell=1}^d \Gamma_{j,\ell}(k)\Gamma_{j',\ell
}(k)=\Gamma_{j,j'}(k),\quad j,j'=1,\ldots,d.
 \end{equation}
The stochastic integral appearing in \eqref{e:10.1}  is a linear
combination of integrals taken with respect to the stochastic
measures, {see
Section \ref{sec3a} below for their precise definitions,
$$
W_{j}^2(ds,dk,ds',dk') :=W_{j}(ds,dk)W_{j}(ds',dk'), \quad
j=1,\ldots,d.
$$}
The structure measure of $W_{j}^2(ds_1,dk_1,ds_2,dk_2)$ equals
\begin{multline}
\bbE\left[W_{j}^2(ds_1,dk_1,ds_2,dk_2)W_{j'}^2(ds_3,dk_3,ds_4,dk_4)\right]=\left[\delta(s_1-s_3)\delta(s_2-s_4)\delta(k_1+ k_3)\delta( k_2+k_4)\vphantom{\int_0^1}\right.\\
\left.\vphantom{\int_0^1}+\delta(s_1-s_4)\delta(s_2-s_3)\delta(k_1+
  k_4)\delta( k_2+k_3)\right]\delta_{j,j'} d{\bf s}_{1,4}d{\bf k}_{1,4}.
  \label{Yab}
\end{multline}
Here we use the notation  $d{\bf s}_{m,n}:=ds_m\ldots ds_n$ and $d{\bf
  k}_{m,n}:=dk_m\ldots dk_n$ for any $m\le n$. 

By \eqref{Yab} and \eqref{e:Gamma2} 
the structure measure of $Y_{j}(ds,dk,ds',dk')$ equals 
\begin{multline}
\label{Ya}
\bbE\left[Y_j(ds_1,dk_1,ds_2,dk_2)Y_{j'}(ds_3,dk_3,ds_4,dk_4)\right]=\left[\delta(s_1-s_3)\delta(s_2-s_4)\delta(k_1+ k_3)\delta( k_2+k_4)\vphantom{\int_0^1}\right.\\
\left.\vphantom{\int_0^1}+\delta(s_1-s_4)\delta(s_2-s_3)\delta(k_1+
  k_4)\delta( k_2+k_3)\right]\Gamma_{j,j'}(k_1+k_2) d{\bf s}_{1,4}d{\bf k}_{1,4}.
\end{multline}
 Note that
the random vector field $V(t,x)$ defined by \eqref{e:10.1} is real
valued and its covariance matrix 
\begin{equation}
R_{j,j'}(t,x)=\bbE[V_{j}(t+t',x+x')V_{j'}(t',x')], \quad
(t,x),(t',x')\in \bbR^{1+d},\,j,j'=1,\ldots,d
\label{e:2.10a}
\end{equation}
is of the form \eqref{010903}, with
\begin{equation}
\label{e:2.10b}
\hat{\cal E}(t,|k|):= |k|^{d-1}\int_{\bbR^d}e^{-[r(|k-\ell|)+r(|\ell|)]|t|}\frac{{a}(|k-\ell|){a}(|\ell|) d\ell}{|k-\ell|^{d+\alpha -1}|\ell|^{d+\alpha-1}},\quad (t,k)\in\bbR^{1+d}.
\end{equation}
\begin{remark}
 \label{rem:1a} { In light of  the assumptions \eqref{Erb} --
\eqref{022806a} we have
$$
\hat{\cal E}(t, \xi)\approx  e^{-\xi^{2\beta}|t|}\frac{a(\xi)}{\xi^{2\al-1}},\quad (t,\xi)\in\bbR\times\R_+.
$$
 Hence the energy spectrum is integrable and the field given by
\eqref{e:10.1} is well defined, provided that $\al<1$.
The diffusivity, defined by the Green-Kubo formula \eqref{TK},
becomes infinite if  $\al+\beta>1$. }
\end{remark}

From \eqref{e:10.1} it is clear that the random field $V$ is stationary and isotropic. Moreover, due to the presence of the term $\Gamma_{j,j'}$ in \eqref{Y} it is also incompressible. 


%

\subsection{Statement of the main result}

Suppose that  $V$ is given by \eqref{e:10.1} and  $x(t)$ is the solution of \eqref{011706a} with $x_0=0$. Let 
  \begin{equation}
   \label{e:H}
   \delta=\frac{\beta}{\al+2\beta-1}\quad \mbox{and}\quad H:=\frac{1}{2\delta}.
  \end{equation}
\begin{remark}
 \label{rem:0} { Note that if $\al+\beta>1$ then $H>1/2$. On the other hand, since
$\al<1$ we have $H<1$.}
\end{remark}

\medskip

Denote
 \begin{equation}x_\varepsilon(t):=x\left(\frac{t}{\varepsilon^{2\delta}}\right)=\varepsilon^{1-2\delta}\int_0^tV\left(\frac{s}{\varepsilon^{2\delta}}, x_\varepsilon(s)\right)ds
\label{e:10.23a}
 \end{equation}
and
\begin{equation}
\label{Einfty}
E_\infty(s,s',k,k'):=a_0r_0\left(\frac{1}{|k||k'|}\right)^{(\al+d-1)/2-\beta}
\exp\left\{-\frac{r_0}{2}\left[|k|^{2\beta}|s|+|k'|^{2\beta}|s'|\right]\right\}.
\end{equation}
Our main result is the following.
\begin{thm}
\label{thm:main}
Suppose that the velocity field $V(\cdot)$ satisfies the assumptions of Section \ref{sec:velocity_field}. Moreover, assume 
 \begin{equation}
 \al<1\qquad \text{and} \qquad \al+\beta>1.
 \label{e:cond_d}
 \end{equation} 
   Then, the processes $\left(x_\varepsilon(t)\right)_{t\ge 0}$ converge in law over $C([0,+\infty);\bbR^d)$, as $\eps\to0$, to the process   $(X(t))_{t\ge 0}=(X_1(t),\ldots, X_d(t))_{t\ge 0}$, given by
\begin{align}
 X_{j}(t):=&\int_0^td\si\int_{(-\infty,\si]^2}\int_{\bbR^{2d}}
  E_\infty(\si-s,\si-s',k,k')
Y_{j}(ds,dk,ds',dk'), \quad j=1,\ldots,d,
\label{e:10.17}
\end{align}
where $Y_j$ is defined in \eqref{Y}.
\end{thm}
The proof of the theorem is presented in Section \ref{sec3}.

\begin{remark}
 \label{rem:1}
{ It is not difficult to see that the process $\left(X(t)\right)$
  is self-similar with index $H$, i.e. the laws of
  $\left(X(ct)\right)_{t\ge 0}$ and $\left(c^HX(t))\right)_{t\ge 0}$ coincide for any
  $c>0$. Moreover, it has stationary increments, that is,
  its law coincides with that of   $\left(X(t+h)-X(h)\right)_{t\ge 0}$ for any
  $h\ge 0$. The process is not Gaussian and takes values in the second Wiener chaos corresponding to the time space white noise $\left(W(dt,dk)\right)$. 
 }
\end{remark}
\begin{remark}
 \label{rem:a1a}
 By stationarity of increments, self-similarity of $X$ and hypercontractivity for double integrals with respect space time white noise  (cf. Proposition \ref{lem:A4} below) for any $m\in \N$ we have
 \begin{equation*}
  \E \abs{X(t+h)-X(t)}^{2m}\le C_m \abs{h}^{2Hm},
 \end{equation*}
hence the process $X$ has H\"older continuous trajectiories with any exponent $\tilde H<H$.
\end{remark}


%

\medskip
 Let us discuss  briefly the relation of the process $\left(X(t)\right)$ with the Rosenblatt process.
Recall (see e.g. \cite{tudor}) that 
the  one dimensional {\em Rosenblatt process} of index $H\in(\frac 12,1)$ is defined as a process of the form 
\begin{equation}
\tilde Z(t):=c\int_{\bbR^2}\left\{\int_0^t(s-y_1)_+^{-(1-H/2)}(s-y_2)_+^{-(1-H/2)}ds\right\}W(dy_1)W(dy_2),\quad t\ge0,
\label{e:Rosenblatt}
\end{equation}
where $W(dy)$ is a one dimensional real valued white noise,
$x_+:=\max(0,x)$, $x\in\bbR$. The process is $H$-self-similar, with
stationary increments.

The $d$-dimensional Rosenblatt process is defined as
$\tilde Z(t):=\left(\tilde Z_1(t),\ldots,\tilde Z_d(t)\right)$, $t\ge0$,  where $\tilde Z_1(\cdot),
\ldots, \tilde Z_d(\cdot)$  are independent one dimensional Rosenblatt processes. 
Let us now consider the process analogous to \eqref{e:10.17} but with $Y_j(ds,dk,ds',dk')$  replaced by $W_j(ds,dk)W_j(ds'dk')$, which corresponds to taking the Kronecker delta $\delta_{j,j'}$ instead of  $\Gamma_{j,j'}(k)$ in \eqref{Y}, that is, let
 \begin{align}
 Z_{j}(t):=&\int_0^td\si\int_{(-\infty,\si]^2}\int_{\bbR^{2d}}
  E_\infty(\si-s,\si-s',k,k')
W_{j}(ds,dk)W_{j}(ds',dk'), \quad j=1,\ldots,d.
\label{e:10.17z}
\end{align}
\begin{proposition}
\label{prop:Rosenblatt}
  The process $Z(t):=\left(Z_1(t),\ldots,Z_d(t)\right)$, $t\ge0$,
  where $Z_j(\cdot)$, $j=1,\ldots,d$, are given by \eqref{e:10.17z}, is a
  $d$-dimensional Rosenblatt process. 
  \end{proposition}
 The proof of the proposition is given in Section \ref{sec5} below. 


\begin{remark}\label{rem:2}
The process $(X(t))$ that appears in the
    statement of Theorem \ref{thm:main} is subordinated to  the
    one defined in \eqref{e:10.17z} in the following sense. From
    representation \eqref{e:10.17z} one concludes that
$Z(t)=X(t)+\tilde X(t)$, where $\tilde X(t)=(\tilde
X_1(t),\ldots,\tilde X_d(t))$ and 
\begin{align}
\tilde X_{j}(t):=\sum_{j'=1}^d&\int_0^td\si\int_{(-\infty,\si]^2}\int_{\bbR^{2d}}
  E_\infty(\si-s,\si-s',k,k')\nonumber\\
&
\times 
|k+k'|^{-2}(k+k')_j(k+k')_{j'}W_{j'}(ds,dk)W_{j'}(ds',dk'), \quad j=1,\ldots,d.
\label{e:10.17b}
\end{align}
The process $\left(\tilde X(t)\right)$
is also self-similar with index $H$ and has stationary increments.

Moreover $\left(X(t)\right)$ and $\left(\tilde X(t)\right)$ are
uncorrelated, i.e.
$\bbE\left [X_j(t)\tilde X_{j'}(t')\right ]=0$ for all $t,t'\ge0$ and
$j,j'=1,\ldots,d$.
\end{remark}

 \section{An interlude on multiple stochastic integrals}

\label{sec3a}

In this section we recall the notion of multiple stochastic integrals
used in the definition of the field $V(\cdot,\cdot)$ and process
$X(\cdot)$. We also discuss some of their properties,
which will be used in the proofs of Theorem \ref{thm:main} and Proposition \ref{prop:Rosenblatt}. We start by describing Feynman diagrams, which provide a  useful tool in our subsequent analysis. 


{Given positive integers $n,r$ we let $Z_{n,r}$ be the set of all pairs
$\{(\ell,m): \ell =1,\ldots, n, m=1,\ldots r\}$. The first element of
a pair $(\ell,m)$  will be called a {\em node}, while the second  one
a {\em hand}.
We  refer to the pair as the $m$-th hand of node $\ell$}. We consider diagrams in which each hand can be either free (not linked to any other), or linked to exactly one hand of a different node (for $r=2$ imagine a group of $n$ people holding hands or not).
We will write $\left((\ell,m),(\ell',m')\right)$ to denote a
link between hand $m$ of node $\ell$ and hand $m'$ of node
$\ell'$. We always require that $\ell\neq \ell'$ and each $(\ell,m)$
can belong to at most one link. The elements of the link are not
ordered but in some calculations it is useful to order them in such a
way that $\ell<\ell'$. Such diagrams will be called {\emph{Feynman}
  diagrams}. They can be equivalently represented as a graph {with
$Z_{n,r}$ being the set of vertices} and the  edges corresponding to the links.
The {vertices with hands  that are not paired} will be called {\emph{free}}.  
In what follows we will represent diagrams described above as a set of links and free vertices. They will be usually denoted by the letter $G$.
The set of all possible diagrams formed over  $n$ nodes with $r$ hands
will be denoted by $\D_n^r$. For $G\in \D_n^r$, we denote by
$G_{free}$, $G_{links}$ the sets of free elements and links,
respectively, while  $f(G)$, $\ell(G)$ stand for 
the cardinalities of the corresponding sets.  Of special importance
are those diagrams, for which $G_{free}=\emptyset$. {We call them {\em complete}}. Following \cite{Talarczyk_silt1} we  denote
them by $\G_n^r:=\{G\in \D_n^r: G_{free} =\emptyset\}$.  Clearly, if $nr$ is odd, then $\G_n^r=\emptyset.$

{  Suppose that  $G\in \G_n^r$. We wish to define a notion of a
  {\em connected component} of $G$, that is a maximal set consisting of those
  links that are connected via common nodes. More precisely,
we say that two nodes $\ell,\ell'$ are 
{\em directly connected} if there exists a link $((\ell,m),(\ell',m'))$, for some $m,m'$.
We say that nodes $\ell$ and $\ell'$ are {\em connected}, and write
$\ell \sim \ell'$ if there exists {a positive integer $k$  and a sequence of nodes $\ell_1,\ldots, \ell_k$ such that $\ell=\ell_1$, $\ell'=\ell_k$ and the node $\ell_{i-1}$ is directly connected to $\ell_{i}$ for each $i=2,\ldots,k$.}
%
%
 It is clear that $\sim$ is an equivalence relation on the set of
 nodes. We say that  $G'\subset G$ is a {\em connected component} of
 $G$ if the nodes of its vertices form an equivalence class of the
 relation $\sim$.  By $\tilde{\G}_n^r$ we denote the subclass of $\G_n^r$ made of those complete Feynman diagrams that have only one connected component.}

Feynman
diagrams $\D_n^1$  can be used to describe Wick products. Diagrams
of the class $\G_n^r$ are useful in computing moments of  Wick products, or of
multiple stochastic integrals (see e.g. \cite{Glimm_Jaffe}, \cite{fox-taqqu},
\cite{Talarczyk_silt1} or formula (3.6) in \cite{janson}).

For example consider the case when $r=1$.
We omit then writing hands in our notation of vertices. Let $f_1,\ldots, f_n\in L^2(\Rdp)$. Suppose that
$(W_1(ds,dk),\ldots.W_d(ds,dk))$ is a white noise as introduced in
Section \ref{sec2.2}. Denote by ${\cal J}_{n,d}^1$ the set of all
multi-indices ${\bf j}:=(j_1,\ldots,j_n)\in \bbZ_d^n$. For each ${\bf j}\in {\cal
  J}_{n,d}^1$ and $f:=f_1\otimes\ldots\otimes f_n$, with $f_j\in
L^2(\bbR^{1+d})$, $j=1,\ldots,n$  we define the multiple stochastic integral
\begin{align}
I_n^{\bf j}&(f)
=
 \int_{\bbR^{(1+d)n}}f_1(s_{1},k_{1})\ldots f_n(s_{n},k_{n})W_{j_1}(ds_{1},dk_{1})\ldots W_{j_n}(ds_{n},dk_{n})\notag\\
& :=\sum_{G\in\D_n^1}(-1)^{\ell(G)}
 \prod_{\left(\ell,\ell'\right)\in G_{links}}\delta_{j_\ell,j_{\ell'}}\left(\int_{\bbR^{1+d}}f_{\ell}(s,k)f_{\ell'}(s,-k)dsdk\right) \prod_{\ell\in G_{free}}\int_{\Rdp} f_{\ell}(s,k)W_{j_\ell}(ds,dk).
\label{e:A2}
 \end{align}
Note that the right hand side of \eqref{e:A2} is simply the Wick
product of the stochastic integrals $\int_{\bbR^{1+d}} f_1dW_{j_1}$, \ldots,
$\int_{\bbR^{1+d}} f_n dW_{j_n}$.
Then 
\begin{equation}
\E\abs{I_n^{\bf j} (f)}^2\le n!\int_{\bbR^{n(1+d)}}|f(s_1,k_1,\ldots,s_n,k_n)|^2 ds_1dk_1\ldots ds_n dk_n.
\label{e:A3}
\end{equation}
The definition of $I_n^{\bf j}$ is extended by linearity to $L^2(\bbR^{n(1+d)})$.
We always have $\E I_n^{\bf j}(f)=0$, $n\ge1$, ${\bf j}\in {\cal
  J}_{n,d}^1$. Moreover $I_n^{\bf j}(f)$ is invariant with respect to permutations of those arguments of $f$ that correspond to the same noises $W_j(\cdot)$.

{Next,} we wish to formulate an analogue of \eqref{e:A2}
for  products of double stochastic integrals.
 This is achieved by an
application of  Proposition 1.1.3, p. 12
of \cite{Nualart}.
Denote by ${\cal J}_{n,d}^2$ the set of all
multi-indices ${\bf j}:=(j_{\ell,m})$, where
{$j_{\ell,m}\in\bbZ_d$, $(\ell,m)\in Z_{n,2}$}.
Suppose that $f\in
L^2(\bbR^{2n(1+d)})$ and
$G\in \D_n^2$. For ${\bf j}=(j_{\ell,m})\in {\cal J}_{n,d}^2$ we define
\begin{align}
 & I_{G,{\bf j}}(f) :=\int_{\bbR^{(1+d)2n}} f(s_{1,1},k_{1,1},s_{1,2},k_{1,2},\ldots, s_{n,1},k_{n,1},s_{n,2},k_{n,2})
  \notag\\
 &\prod_{\left((\ell,m),(\ell',m')\right)\in G_{links}} \delta_{j_{\ell,m},j_{\ell',m'}}\delta(s_{\ell,m}-s_{\ell',m'})\delta(k_{\ell,m}+k_{\ell',m'})ds_{\ell,m}ds_{\ell',m'}dk_{\ell,m}dk_{\ell',m'}\notag\\
 &\hskip 0.6cm \prod_{(\ell,m)\in G_{free}}W_{j_{\ell,m}}(ds_{\ell,m},dk_{\ell,m}).
 \label{e:A5}
\end{align}
The above formula can be interpreted as follows: to compute  $I_{G,\bf
  j}(f)$   we identify the variables corresponding to a link  {$\left((\ell,m),(\ell',m')\right)$ and
integrate out according to \eqref{010403}. The remaining variables,
corresponding to free vertices $(\ell,m)$ in $G$,} are integrated with respect to the
corresponding noises $W_{j_{\ell,m}}(\cdot)$. The resulting   multiple
stochastic integral
of order $f(G)$ is interpreted via formula \eqref{e:A2}.
The following   result  
is a direct consequence of Proposition 1.1.3, p. 12
of \cite{Nualart}.
\begin{proposition}
\label{lem:A1}
 {Suppose that $f_1,\ldots, f_n\in L^2(\bbR^{2(1+d)})$. Let ${\bf j}=(j_{\ell,m})$
 be a multi-index, as described in the foregoing.} Then
 \begin{equation}\label{e:A6}
 I_2^{j_{1,1},j_{1,2}}(f_1)\ldots I_2^{j_{n,1},j_{n,2}}(f_n)=\sum_{G\in \D_n^2}I_{G,{\bf j}}(f_1\otimes \ldots \otimes f_n).
 \end{equation}
\end{proposition}

\medskip
Using the fact that the mean of any multiple stochastic integral vanishes, 
one immediately obtains, as a corollary from Proposition \ref{lem:A1},
a well known formula for {  the moment of the (Wick) product  of
  $n$ double integrals} (see e.g. Lemma 2.1 in \cite{Talarczyk_silt1}).
\begin{corollary}
 \label{lem:A2} Under the assumptions of Proposition $\ref{lem:A1}$ we have
 \begin{equation}
 \E\left[\vphantom{\int_0^1} I_2^{j_{1,1},j_{1,2}}(f_1)\ldots
   I_2^{j_{n,1},j_{n,2}}(f_n)\right]=\sum_{G\in \G_n^2}I_{G,{\bf j}}(f_1\otimes \ldots \otimes f_n).
 \label{e:A7}
 \end{equation}
 \end{corollary}

 The moments of multiple itegrals can be estimated with the help of the
 second moment, thanks to Theorem 3.50, p. 39 of \cite{janson}:
\begin{proposition}\label{lem:A4}
  For any $p\ge 1$ there exists $C_p>0$ such that for any $f\in
  L^2(\R^{(d+1)n})$ and multi-index ${\bf j}\in{\cal J}_{n,d}^2$ we have
 \begin{equation}\label{e:A9}
\left(  \E\abs{I_n^{{\bf j}}(f)}^p\right)^\frac 1p\le C_p \left(\E\abs{ I_n^{{\bf j}}(f)}^2\right)^\frac 12.
 \end{equation}
\end{proposition}

Finally, we formulate a result on conditioning, which is derived in exactly the same way as Lemma 1.2.5 in \cite{Nualart}.  Denote 
\begin{equation}
\label{Fs}
\F_s=\sigma\{\langle \phi,W_j\rangle,\quad \phi\in
L^2(\bbR^{1+d}),\,{\rm supp }\,\phi\in (-\infty,s]\times\bbR^d,\,j=1,\ldots,d\}.
\end{equation}
\begin{proposition}
 \label{lem:A3}
and let $f\in L^2(\R^{(d+1)n})$.
Then, 
\begin{equation}
 \label{e:A8}
 \E(I_n^{{\bf j}}(f)|\F_s)=I_n^{{\bf j}}(f\ind_{(-\infty,s]\times
   \Rd}^{\otimes n}),\quad \mbox{ for any multi-index }{\bf j}\in {\cal J}_{n,d}^2.
\end{equation}
\end{proposition}

 \section{Proof of Theorem \ref{thm:main}}

\label{sec3}
\subsection{{Outline} of the proof}In this section we describe the main steps of the proof. The proofs of some technical lemmas are presented in the following sections.
\subsubsection{Reformulation of the problem}

The first step is to reformulate our problem. 
Let $T:=\varepsilon^{-\delta/\beta}$ and  $\tilde x_T(t):=x_\varepsilon(t)$. By \eqref{e:H} and \eqref{e:10.23a} we see that  $\tilde x_T$ satisfies
 \begin{equation}
  \tilde x_T(t)=T^{2\beta (1-H)}\int_0^{t}V(sT^{2\beta},\tilde x_T(s))ds.
\label{e:10.23b}
  \end{equation}
 {Our goal is} to prove the convergence in law of
 $\tilde x_T(\cdot)$, as $T\to \infty$. 

\bigskip


\bigskip

Let $V_T=(V_{T,1},\ldots, V_{T,d})$, where
\begin{align}
\label{VT}
V_{T,j}(t,x):=&\sum_{j'=1}^d V_{T,j,j'}(t,x)\quad \mbox{and}\\
V_{T,j,j'}(t,x):=&\int_{\bbR^{2d}}\int_{(-\infty,t]^2}
 \exp\left\{\ii(k+k')\cdot
  x\right\} E_T(t-s,t-s',k,k')\notag\\
  &\hskip 2cm\times
\Gamma_{j,j'}(k+k') W_{j'}(ds,dk) W_{j'}(ds',dk'),
\quad (t,x)\in\bbR^{1+d},\label{e:VTj},
\end{align}
with
\begin{equation}
\label{ET}
E_T(s,s',k,k'):=
\left(\frac {{ r_T(\abs{k})a(\abs{ k/T})r_T(\abs{k'})a(\abs{ {k'}/T}) }}
{|k'|^{(d+\alpha -1)}|k|^{(d+\alpha-1)}}\right)^{\frac 12}
\exp\left\{-\frac12\left[r_T(|k|)|s|+r_T(|k'|)|s'|\right]\right\},
\end{equation} 
\begin{equation}\label{rT}
r_T(\xi):=T^{2\beta}r(T^{-1}\xi).
\end{equation}
Recalling \eqref{Erb}--\eqref{e:Gamma} and \eqref{Einfty}, and using the fact that
\begin{equation}
W(T^{2\beta}ds,T^{-1}dk)
\overset{d}{=}T^{2\beta/2-d/2}W(ds,dk),
\label{e:WT_scaling}
\end{equation} 
where $\overset{d}{=}$
denotes the equality of the laws of random fields, we can easily
conclude the following lemma.
\begin{lemma}
 \label{lem:10.2} Assume that 
  $V(\cdot)$ is the random field defined by \eqref{e:10.1} and {satisfies} the assumptions {made} in Section $\ref{sec:velocity_field}$. Let $H$ be given by \eqref{e:H}.
Then for any $T>0$ 
\begin{equation}
\label{e:10.9}
\left( T^{2\beta(1-H)}V(T^{2\beta}t, Tx)\right)_{(t,x)\in\bbR^{1+d}}\overset{d}{=}(V_T(t,x))_{(t,x)\in\bbR^{1+d}}.
\end{equation}
\end{lemma}
The proof of this lemma is straightforward, using \eqref{e:WT_scaling}, therefore we omit it.

 Directly from \eqref{e:10.23b} and Lemma \ref{lem:10.2} we
 {obtain} the following.
 \begin{corollary}\label{lem:10.9}  Under the assumptions of Theorem $\ref{thm:main}$ and for $T=\varepsilon^{-\delta/\beta}$ the process $x_{\varepsilon}(\cdot)$ has the same law as $z_T(\cdot)$, where
 \begin{equation}
 z_T(t)=\int_0^t V_T\left(s, \frac{z_T(s)}{T}\right)ds, \quad t\in\bbR,
 \label{e:10.24}
 \end{equation}
with $V_T$ defined by \eqref{VT}--\eqref{rT}.
\end{corollary}
Hence, to prove Theorem \ref{thm:main} it suffices to show that, as
$T\to \infty$, the processes $z_T(\cdot)$ converge in law in
$C([0,\infty),\Rd)$ to  the process $X(\cdot)$ defined in
\eqref{e:10.17}. The general approach is standard: we prove the
convergence of finite dimensional distributions {(in fact we show the
moment convergence of finite dimensional distributions)}  and then
{proceed with establishing} tightness.

\subsubsection{A result on stationarity}
Note that $V_T(\cdot)$  given by \eqref{VT}-\eqref{rT} is again
stationary and divergence free field. It is {well}
known that in this case the process  $\left(V_T(s,z_T(s)/T)\right)_{s\in \R}$,
{with $z_T(\cdot)$ given by \eqref{e:10.24}}, is
stationary (see \cite{port-stone}). In the course of the proof we will
need a {somewhat stronger} result

Given  $p\in \bbZ_+$ consider 
the family of fields of the form
\begin{equation}
{\cal D}^p(V_T):=\left\{\frac{1}{T^p}\partial_{z_1}^{p_1}\ldots\partial_{z_d}^{p_d}V_{T,j},\,j=1,\ldots,d,\,p_1,\ldots, p_d\in\bbZ_+,\sum_{k=1}^dp_k=p\right\}.
\label{e:DpVT}
\end{equation}
We allow $p=0$ with the convention ${\cal D}^0(V_T)=\{V_{T,j}, j=1,\ldots, d\}$.
 Since the spatial realizations of $V_T(t,x)$ are incompressible,  a direct application of Theorem
2, p. 500 of \cite{port-stone} yields the following.
\begin{proposition}
\label{prop011007} Let $z_T$ be given by \eqref{e:10.24}.
For any $T>0$, $p_j\in \bbZ_+$, $s_j,s_j'\in \R$, $H_j\in {\cal D}^{p_j}(V_T)$,  $j=1,\ldots,m$ the
process
$$
\left( \vphantom{\int_0^1}H_1\left(s_1+t,z_T(s_1'+t)/T\right),\ldots, H_m\left(s_m+t,z_T(s_m'+t)/T\right)\right)_{t\in\bbR}
$$
is stationary.
\end{proposition}

\subsubsection{Convergence of finite dimensional distributions}
\medskip
Thanks to the reformulation given in Corollary \ref{lem:10.9} it is
quite easy to guess what the limit should be. {Since $z_T(s)/T$ is
expected to become small, as $T\to+\infty$, formula \eqref{e:10.24} suggests that
$z_T(\cdot)$ should be close to the process }
\begin{equation}
y_T(t):=\int_0^t V_T(s,0)ds, \qquad t\ge 0.\label{e:yT}
\end{equation}
Note that by \eqref{rT} and  \eqref{022806a}  {we have}
\begin{equation}
\label{e:crT}\lim_{T\to \infty}r_T(\xi)=r_0\abs{\xi}^{2\beta}\qquad \textrm{and } \qquad 
r_T(\xi)\le C\abs{\xi}^{2\beta},\qquad \xi\in \Rd.
\end{equation}
Hence, { using  the continuity of $a(\cdot)$ at $0$ and \eqref{022806a}}, we conclude
\begin{equation}
\label{e:limET}
\lim_{T\to+\infty}E_T(s,s',k,k')=E_\infty(s,s',k,k')
\end{equation}
 pointwise  (cf \eqref{ET} and \eqref{Einfty}). Thanks to \eqref{Erb}, boundedness of $a$ and \eqref{e:crT}, we obtain
 \begin{equation}
\label{010707}
E_T(s,s',k,k')\preceq E_\infty(s,s',k,k'),\quad (s,s',k,k')\in \bar\bbR^2_+\times\bbR^{2d},\,T\ge1.
\end{equation}
Hence $y_T(\cdot)$ should converge to $X(\cdot)$, defined in
\eqref{e:10.17}, which indeed  is  the case.
\begin{proposition}
 \label{lem:10.6} Let $y_T$ be defined by \eqref{e:yT} and \eqref{VT}--\eqref{rT}. Then under the assumptions of Theorem \ref{thm:main} for any $q>0$ we have
  \begin{equation}
  \lim_{T\to +\infty}\E \abs{y_T(t)-X(t)}^{q}=0,\quad t>0.
  \label{e:11.1}
 \end{equation}
\end{proposition}
The proof of the above proposition is presented in  Section
\ref{sec:4.2}.

 {The argument presented in the foregoing justifies
heuristically the validity of the claim made in Theorem 
\ref{thm:main}.}
Note, however, that the field $V_T(\cdot)$ itself does not converge,
as $T\to+\infty$. We have the following lemma, that will be useful
later {on.}
\begin{lemma}\label{lem:10.9aa} Suppose that the assumptions of
  Theorem \ref{thm:main} are satisfied. { Then, the following are true.}\\
a) There exists $C>0$ such that
  \begin{equation}
\bbE \abs{V_T(0,0)}^2=\sum_{j=1}^d\bbE V_{T,j}^2\left(0,0\right)=C
T^{2(1-\alpha)},\quad T>0.
 \label{e:10.24aa}
 \end{equation}
 b) For any $n\in \N$,  $q>0$ and $j, j_1\ldots , j_n\in \{1,\ldots, d\}$ there exists $C_q>0$ such that 
 \begin{equation}
\left(\bbE\abs{ \frac 1{T^n}\partial^n_{{x_{j_1},\ldots, x_{j_n}}} V_{T,j}\left(0,0\right)}^q\right)^{\frac 1q}\le C_q T^{(1-\alpha)},\quad \mbox{as }T{\ge}1.
 \label{e:10.24aaa}
 \end{equation}
 c) For any $n\in \N$,  $j'\in \{1,\ldots, d\}$ there exists $C>0$ such that
 \begin{equation}\E \left( \frac 1{T}\int_0^t \partial_{x_{j'}} V_{T,j}(s,0)ds\right )^2
 \le Ct(t\vee 1) \begin{cases}
 T^{-2(\alpha+\beta-1)} ,\qquad &\textrm{if} \quad \alpha+\beta<2,\\
 T^{-2}\log T, &\textrm{if }\quad \alpha+\beta=2,\\
 T^{-2},\quad &\textrm {if} \quad\alpha+\beta>2.
 \end{cases}
 \label{e:deriv}
 \end{equation}
\end{lemma}
The proof of this lemma is given in Section \ref{sec:4.3}

\begin{remark}
 Using Lemma \ref{lem:10.9aa} we can easily conclude the convergence
 of finite dimensional distributions in Theorem \ref{thm:main} for a
 narrower range of parameters { when $\alpha<1$ and $\alpha+ \beta/2>1$.
 By \eqref{e:10.24} and \eqref{e:yT}
 we can write
 \begin{equation}
  \label{e:zy}
    z_{T,\ell}(t)-y_{T,\ell}(t)=\sum_{j=1}^d
\int_0^tds\int_0^s \frac 1T \partial_{x_j}{V_{T,\ell}}\left(s,\frac{z_T(u)}T\right)
 V_{T,j}\left (u,\frac{z_T(u)}T\right) du .
 \end{equation}
Changing the order of integration and  applying Cauchy-Schwarz inequality we have
 \begin{equation*}
  \E \left[\abs{z_{T,\ell}(t)-y_{T,\ell}(t)}\right]
 \le \sum_{j=1}^d \left\{\E \int_0^tdu\left[\int_u^t \frac 1T \partial_{x_j}{V_{T,\ell}}\left(s,\frac{z_T(u)}T\right)ds \right]^2
   \E\int_0^t V_{T,j}^2 \left (u,\frac{z_T(u)}T) \right)du\right\}^{\frac 12}.
 \end{equation*}
 By Proposition \ref{prop011007} 
 the processes $\left(V_T(u,z_T(u)/T)\right)_{u\ge 0}$ and
 $\left(\partial_{x_j}{V_{T,\ell}}(r+u,z_T(u)/T)\right)_{u\ge 0}$ are stationary, hence
 \begin{equation}
 \label{e:14.1}
  \E \left[\abs{z_T(t)-y_T(t)}\right]\le  \sqrt  t \sum_{j,\ell=1}^d \left\{\E \int_0^t\left[\int_{0}^{t-u} \frac1T\partial_{x_j}V_{T,\ell}(s,0) ds\right]^2du\right\}^{\frac 12} 
\left\{\vphantom{\int_)^1}\E V_{T,j}^2(0,0 )\right\}^{\frac 12}.
 \end{equation}
By parts a) and c) of  Lemma \ref{lem:10.9aa} we conclude that  the right hand
side of \eqref{e:14.1} converges to $0$, as $T\to \infty$,  provided
that $\alpha+{\beta}/{2}>1$ and $\al<1$. This together with
\eqref{e:11.1} yields in particular
the weak convergence of finite dimensional distributions  of $z_T(\cdot)$
to $X(\cdot)$, in this case. By hypercontractivity (see Proposition
\ref{lem:A4}) one can also argue that the finite dimensional
distributions of $z_T(\cdot)$ converge to those of
$X(\cdot)$, in the $L^q$ sense for any $q>0$.}

{Quite remarkably, the extension of the above argument to
the case $\alpha+\beta>1$, $\al<1$ (considered in Theorem
\ref{thm:main}) eludes us, although we believe that the finite
dimensional distributions of $z_T(\cdot)$ should converge in the $L^q$
sense for any $q>0$ then as well. However,  the natural idea of using
 a longer Taylor expansion and applying the Cauchy-Schwarz (or
 H\"older) inequality to estimate the remainder does not seem to
 work. Therefore, to show the convergence of finite dimensional distributions we use the method of moments. There we also use the Taylor expansion and  the Cauchy-Schwarz inequality but in a more subtle way. }\end{remark}



It is known (see e.g. p. 113 in \cite{fox-taqqu}) that the law of  a random variable belonging to the second Wiener chaos is determined by its moments (this is no longer true in chaoses of higher order).
Therefore, to prove convergence of finite dimensional distributions it suffices to show that for any $m,n\in \N$, $0<t_n\le \ldots \le t_1$ and $a_{r,j}\in \R$, $r=1,\ldots, n$, $j=1,\ldots, d$ we have
\begin{equation}
 \lim_{T\to \infty}\E\left( \sum_{r=1}^n\sum_{j=1}^d a_{r,j}z_{T,j}(t_r)\right)^m
 =\E\left( \sum_{r=1}^n\sum_{j=1}^d a_{r,j}X_{j}(t_r)\right)^m.
 \label{e:14.2}
\end{equation}
{The latter is a  direct consequence of the following.}

%
%
%
%
%
\begin{thm}[Convergence of moments]
\label{thm090801}
Suppose that  the assumptions of Theorem $\ref{thm:main}$ are satisfied and $z_T$ is defined by \eqref{e:10.24}. Then for any $t_1\ge t_2\ge\ldots\ge t_m\ge0$,
 non-negative integers 
$j_1,\ldots,j_m\in\bbZ_d$ we have
\begin{align}
\label{010306x}
\lim_{T\to+\infty}\bbE\left[ \prod_{n=1}^mz_{T,j_n}(t_n)\right]=\bbE\left[ \prod_{n=1}^mX_{j_n}(t_n)\right].
\end{align}
\label{thm011206-18y}
\end{thm}
 {The proof of the above result is the most difficult
   part of the argument. Its main steps are presented  in Section \ref{sec:proof_moments},
 leaving  the demonstration of some technical lemmas till Sections \ref{sec5.5-18} and \ref{sec5.4-18}.
 The key ingredients of the proof of Theorem \ref{thm090801} are: inductive Taylor expansions (using
 \eqref{e:10.24}), applications of condtitioning and stationarity, resulting
 from the incompressibility of $V(\cdot)$ (following from Proposition
 \ref{prop011007}) and   Propositions \ref{lem:A1}--\ref{lem:A3}
 applied to calculate the expectations (both conditional and
 unconditional) of the expressions resulting from
 the Taylor expansion}.
 
 As mentioned above, from Theorem \ref{thm011206-18y} we obtain {in particular}:
 \begin{corollary}[Convergence of finite dimensional distributions]
For any $t_1\ge t_2\ge\ldots\ge t_m\ge0$ the laws of  $\left(z_{T}(t_1),\ldots,z_{T}(t_m)\right)$ converge weakly to the law of 
$\left(X(t_1),\ldots,X(t_m)\right)$.
\label{cor:fin_dim}
\end{corollary}


\subsubsection{Tightness}

Finally, the last step is to prove tightness. We show 
\begin{proposition}
 \label{prop:tightness}
 Under the assumptions of Theorem $\ref{thm:main}$
 there exists $C>0$ such that 
 \begin{equation*}
  \E \left|z_{T}(s+t)-z_{T}(s)\right|^2\le Ct^{2H},\quad\mbox{ for any }s,t\ge 0,\,T>1.
 \end{equation*}
\end{proposition}
{The above results implies tightness of the family of
  laws of $z_T(\cdot)$, $T>1$ in $C([0,+\infty);\bbR^d)$, provided
  that $H>1/2$, see e.g. Theorem 12.3, p. 95 of \cite{billingsley}.}
The proof of Proposition \ref{prop:tightness} is presented in Section \ref{sec:tightness}. It uses  techniques developed in the proof of Theorem \ref{thm011206-18y}.

Combining Corollaries \ref{cor:fin_dim} and \ref{lem:10.9} with Proposition \ref{prop:tightness}, recalling that $H>\frac 12$ (see Remark \ref{rem:0}), concludes the proof of Theorem \ref{thm:main}.

\subsection{Proof of Proposition \ref{lem:10.6}}

\label{sec:4.2}

%

Using  the definitions of the
processes  $X(\cdot)$, $z_T(\cdot)$ and the field $V_T(\cdot)$, see \eqref{e:10.17} and \eqref{VT}-\eqref{e:10.24},
and the $L^2$ isometry for multiple integrals we have
\begin{align*}
 &\sum_{j=1}^d\E \left[y_{T,j}(t)-X_j(t)\right]^2 =
 2(d-1)\int_{\bbR^{2d+2}}dkdk'dsds'\\
&\times\left[\int_0^t \ind_{\{s\le r, s'\le r\}}
\left( \vphantom{\int_0^1}E_T(r-s,r-s',k,k')- E_\infty(r-s,r-s',k,k')\right) dr\right]^2 
\end{align*}
(cf \eqref{ET} and \eqref{Einfty}).
By \eqref{e:limET} and \eqref{010707} we can {use} the Lebesgue dominated
convergence theorem to conclude \eqref{e:11.1} for $q=2$, provided we can show that
\begin{equation}
 \label{e:11.3}
 J_c:=\int_{\bbR^{2d+2}}dkdk'dsds'
\left[\int_0^t \ind_{\{s\le r, s'\le r\}}
E_\infty(c(r-s),c(r-s'),k,k')dr\right]^2 <+\infty,\quad c>0.
\end{equation}
Integrating with respect $s$ and $s'$ 
and using
\begin{equation}
 \label{e:integral}
 \int_{-\infty}^{r\wedge r'}
 \abs{k}^{2\beta}e^{-\frac {cr_0}2 (r+r'-2s)\abs{k}^{2\beta}}ds=\frac 1{cr_0}e^{-\frac {cr_0}2 \abs{r-r'}\abs{k}^{2\beta}}
\end{equation}
we conclude that $J_c$ equals,
up to a constant, 
 \begin{equation}
\label{010807}
  \int_{\R^{2d}}\int_{[0,t]^2}\frac{1}{(\abs{k}\abs{k'})^{\al+d-1}}\exp\left\{-\frac{cr_0}{2}\abs{r-r'}\left(\abs{k}^{2\beta}+\abs{k'}^{2\beta}\right)\right\}drdr'dkdk'.
 \end{equation}
We have the following elementary inequality
\begin{equation}
\label{021707}
\frac{1}{\ga}(1-e^{-\gamma t}) (1+\ga)\le C(t):=2 (t\vee 1),\quad \ga,t>0.
\end{equation}
Note that for any $\ga,t>0$ we have
\begin{equation}
 \int_{[0,t]^2} e^{-\gamma
   \abs{r-r'}}drdr'=\frac{2}{\ga}\int_0^t(1-e^{-\ga r})dr\le 2t\frac {1-e^{-\gamma t}}\gamma
\le \frac {2tC(t)}{1+\gamma}.
 \label{e:11.10a}
 \end{equation}
Using an elementary inequality
 \begin{equation}
  \frac 1{1+a^2+b^2}\le \frac 2{(1+a)(1+b)},\qquad a,b>-1
 \label{e:11.10b}
 \end{equation}
we conclude that the expression in \eqref{010807} is bounded from
above, modulo some absolute constant, by
 \begin{equation}
\frak{c}(t)\int_{\R^{2d}}\frac{ dkdk' }{(\abs{k}\abs{k'})^{\al+d-1}(1+\abs{k}^{2\beta}+\abs{k'}^{2\beta})}
  \le \frak{c}(t) \left[\int_\Rd\frac 
{dk}{\abs{k}^{\al+d-1}(1+\abs{k}^\beta)}\right]^2<+\infty,
\label{e:11.11a}
  \end{equation}
due to the second inequality in \eqref{e:cond_d}. Here
\begin{equation}
\label{frak-c}
\frak{c}(t):=t(t\vee 1),\quad t\in\bar\bbR_+.
\end{equation}
This finishes the
proof of \eqref{e:11.3}, thus \eqref{e:11.1} follows for $q=2$. The result for an arbitrary $q\ge1$ then follows
from the fact that $y_{T}(t)-X(t)$ belongs to the second  Gaussian chaos space, corresponding to $\left(W(\cdot)\right)$, where all the $L^q$ norms are equivalent, see Lemma \ref{lem:A4}. \qed

\subsection{Proof of Lemma \ref{lem:10.9aa}}
\label{sec:4.3}

%
%
%
%
%
%
%

a) {In the same way as} in the proof of Proposition \ref{lem:10.6} we have
\begin{equation*}
 \E V_T^2(0,0)=2(d-1)\int_{\R^{2d+2}}\ind_{[s\le 0, s'\le 0]}E_T^2(-s,-s',k,k')dkdk'dsds',
\end{equation*}
with $E_T$ defined by \eqref{ET}.
Integrating with respect to $s$ and $s'$, similarly as in \eqref{e:integral}, we obtain
\begin{equation*}
 \E V_T^2(0,0)=2(d-1)\left(\int_{\Rd}\frac{a(\abs{kT^{-1}})}{\abs{k}^{d+\alpha-1}}dk\right)^2 =2(d-1)T^{2(1-\alpha)}\left(\int_\Rd \frac{a(\abs{k})}{\abs{k}^{d+\alpha-1}}dk\right)^2.
\end{equation*}
{Note that the integral on the right hand side is
  finite, since $\alpha<1$ and $a(\cdot)$ has a compact support.
This ends  the proof of part a) of the Lemma. }

b) It is clear that it suffices to show the estimate for each component $V_{T,j,j'}$ of $V_{t,j}$ (cf. \eqref{VT}-\eqref{e:VTj}). By \eqref{e:VTj} we have
\begin{multline}
 \frac 1{T^n}\partial^n_{x_{j_1}\ldots, x_{j_n}}V_{T,j,j'}(t,0)\\
 ={ \rm i}^n\int_{\R^{2d}}\int_{(-\infty,t]^2}
\left[\prod_{\ell=1}^n{\frac{(k+k') _{j_\ell} }{T}}\right]E_T(t-s,t-s',k,k')
\Gamma_{j,j'}(k+k') W_{j'}(ds,dk) W_{j'}(ds',dk')
\label{e:dVj}
\end{multline}
Note that $\abs{\Gamma_{j,j'}(k+k')}\le 1$ and the function $k\mapsto \abs{k}^na(\abs{k})$ is bounded and  compactly supported.  
Hence, the second moment  can be estimated as in part $a)$ of the
lemma. Thus, {the assertion of part b)} holds for $q=2$. The general case follows by hypercontractivity  (see Proposition \ref{lem:A4}).

c) By stationarity we have
\begin{equation*}
 \E \abs{\frac1T\int_0^t\partial_{x_m} V_{T,j}(s,0)ds}^2={\frac{2}{T^2}\sum_{j=1}^d\int_0^tdu\int_0^u\E\left[\partial_{x_m} V_{T,j}(s,0)\partial_{x_m} V_{T,j}(0,0)\right]ds}.
\end{equation*}
{Applying \eqref{e:dVj} and estimates \eqref{e:A2} and \eqref{e:crT}
we obtain
\begin{align*}
  \E &\abs{\frac1T\int_0^t\partial_{x_m} V_{T,j}(s,0)ds}^2\preceq\int_0^t\int_0^u\int_{\R^{2d}} \abs{kT^{-1}}^2\frac
{a(\abs{kT^{-1}})a(\abs{k'T^{-1}})}
       {\abs{k}^{d+\alpha-1}\abs{k'}^{d+\alpha-1}}e^{-cs
       (\abs{k}^{2\beta}+\abs{k'}^{2\beta})}dkdk'dsdu
\end{align*}
for some constant $c>0$.
Here we have also used the fact that
 $|\Gamma_{j,j'}(k+k')|\le 1$ and an elementary inequality
$$
\abs{\frac{k+k'}{T}}^2\le 2\left(\left|\frac{k}{T}\right|^2+\left|\frac{k'}{T}\right|^2\right)
$$
Thanks to the fact that  $a(\cdot)$ is bounded and  \eqref{021707}
we conclude that
\begin{align*}
& \E \abs{\frac1T\int_0^t\partial_{x_m} V_{T,j}(s,0)ds}^2 \preceq {\frak c}(t)\int_{\R^{2d}}\frac
{\abs{kT^{-1}}^2 a(\abs{kT^{-1}})} {\abs{k}^{d+\alpha-1}\abs{k'}^{d+\alpha-1}(1+\abs{k}^{2\beta}+\abs{k'}^{2\beta})}dkdk'.
\end{align*}
where ${\frak c}(t)$ is given by \eqref{frak-c}.
Substituting $k'\mapsto {k'}{(1+\abs{k}^{2\beta})^{1/{2\beta}}}$ we get
\begin{equation*}
 \int_{\Rd}\frac{ dk'}{\abs{k'}^{d+\alpha-1} (1+\abs{k}^{2\beta}+\abs{k'}^{2\beta})}
 =\frac C{(1+\abs{k}^{2\beta})^{\frac {\alpha-1+2\beta}{2\beta}}
 }
\end{equation*}
for some constant $C>0$.
Summarizing, we have shown so far that 
\begin{equation*}
  \sum_{m=1}^d\E \abs{\frac1T\int_0^t\partial_{x_m} V_{T}(s,0)ds}^2\preceq \frak
  c(t)I(T),\quad T,t>0,
 \end{equation*}
 where
 \begin{equation*}
I(T):=  \int_\Rd\frac{\abs{kT^{-1}}^2a(\abs{kT^{-1}})}{\abs{k}^{d+\alpha-1}(1+\abs{k}^{2\beta})^{\frac {\alpha-1+2\beta}{2\beta}}}dk.
\end{equation*}}
If $\alpha+\beta>2$ we use the fact that  $a(\cdot)$ is bounded and estimate
\begin{equation*}
 I(T)\preceq\frac 1{T^2}\int_\Rd\frac{dk}{\abs{k}^{d+\alpha-3}(1+\abs{k}^{2\beta})^{\frac {\alpha-1+2\beta}{2\beta}}}=C T^{-2}
\end{equation*}
{for some constant $C>0$.}
For $\alpha+\beta<2$, we   use the fact that  $a(\cdot)$ is compactly
supported, therefore
\begin{equation*}
 I(T)\le\frac 1{T^2}\int_\Rd\frac{a(\abs{kT^{-1}})}{\abs{k}^{d+2\alpha+2\beta -4}}dk
 =C T^{-(2\alpha+2\beta-2)}.
\end{equation*}
{Finally, if $\alpha+\beta=2$ then, assuming that ${\rm
  supp}\,a(\cdot)\subset[0,K]$ for some $K>0$, we can write
\begin{align*}
& I(T)\preceq  T^{-2}\left(\int_{\abs{k}\le 1}\frac {dk}{\abs{k}^{d+\alpha-3}}+\int_{1<\abs{k}<KT} \frac{dk}{\abs{k}^d} \right)\\
& \preceq T^{-2}\left( 1+\log T\right)
\preceq CT^{-2}\log T
\end{align*}
for $T\ge  1$.}
This concludes the proof of the lemma.
\qed

\subsection{Proof of Theorem  \ref{thm011206-18y}}

\label{sec:proof_moments}
Recall that  
$z_T(\cdot)$ is defined by  {\eqref{e:10.24}}.
By Proposition \ref{lem:10.6} to prove \eqref{010306x} it suffices to show
%
 that for any $t_1\ge t_2\ge \ldots t_m\ge 0$  and
$j_1,\ldots,j_m\in\bbZ_d$  (cf \eqref{010810}):
\begin{align}
\label{010306}
\lim_{T\to+\infty}\left\{\bbE \int\displaylimits_{\Box({\bf t}_{1,m})}\prod_{n=1}^m V_{T,j_n}\left(s_n, \frac {z_T(s_n)}{T}\right)d{\bf s}_{1,n}-\bbE \int\displaylimits_{\Box({\bf t}_{1,m})}\prod_{n=1}^m V_{T,j_n}\left(s_n, 0\right)d{ \bf s}_{1,n}\right\}=0.
\end{align}

%

If we split  the multiple integrals in  \eqref{010306} depending on the order of the variables $s_n$, and note that 
 for any permutation {$(\pi_n)_{1\le n\le m}$} of the set $\bbZ_m$ the set
\begin{equation*}
 \{{\bf s}_{1,m}\in\R_+^m : s_{\pi_1}\le s_{\pi_2}\le \ldots\le s_{\pi_m}, s_{\pi_n}\le t_{\pi_n},n=1,\ldots m \}
  \end{equation*}
is equal to
 \begin{equation*}
 \{{\bf s}_{1,m}\in\R_+^m : s_{\pi_1}\le s_{\pi_2}\le \ldots\le s_{\pi_m}, s_{\pi_n}\le t_{\pi_n}\wedge t_{\pi_{n+1}}\ldots \wedge t_{\pi_m}, n=1,2\ldots, m\}
\end{equation*} 
then it is clear that to prove \eqref{010306} it suffices to show that
\begin{align}
\label{010306az}
\lim_{T\to+\infty}\left\{\int_{\Delta({\bf t}_{1,m})}\bbE\left[\prod_{n=1}^mV_{T,j_n}\left(s_n, \frac {z_T(s_n)}{T}\right)\right]d{\bf s}_{1,m}-\int_{\Delta({\bf t}_{1,m})}\bbE\left[\prod_{n=1}^m
  V_{T,j_n}\left(s_n, 0\right)\right]d{\bf s}_{1,m}\right\}=0.
\end{align}

Let us consider the first integral in \eqref{010306az}. The idea is
 to
write first {the} Taylor expansion of $V_{T,j_1}(s_1, z_T(s_1)/T)$ in the second variable  at $z_T(s_2)/T$ using the equation  \eqref{e:10.24}.
{It turns out that if the expansion is long enough,
  then  we can   show that  the remainder term converges to $0$.} {To
calculate the limit of the other   terms we perform the Taylor
expansion with respect to $s_2$ at $z_T(s_3)/T$ and repeat the
procedure untill we reach the variable $s_m$, expanding around
$0$. Eventually, discarding the remainder terms we arrive at
expressions without a random spatial argument. After some expilict
calculation we conclude in this way equality \eqref{010306az}.}

%

To make this precise we need several lemmas and some additional notation. 
Given $s_1\ge\ldots\ge s_m$ and $H_j\in {\cal D}^{p_j}(V_T)$, $p_j\in \Z_+$, $j=1,\ldots,m$ (see \eqref{e:DpVT}) let
\begin{equation}
\label{H0}
H^{(0)}\left({\bf s}_{1,m},z\right):=\prod_{j=1}^m
  H_{j}\left(s_j,z\right).
\end{equation}
Recall that
${\bf s}_{1,m}:=(s_{1},\ldots,s_m)$.
Suppose that  $H^{(N)}({\bf s}_{1,m+N},z)$ has been defined for some $N\ge0$. For $s_1\ge \ldots s_{m+N+1}$ we let 
\begin{align}
\label{e:W_n1x}
H^{(N+1)}({\bf s}_{1,m+N+1},z):=\frac 1{T}\sum_{j'=1}^d\partial_{z_{j}}H^{(N)}({\bf s}_{1,m+N},z) V_{T,j}(s_{m+N+1},z).
 \end{align}
Suppose that  $s_1\ge \ldots\ge s_m\ge s'$. Using the above notation  we can write the Taylor expansion of $H^{(0)}(s_{1,m},z_T(s_m)/T)$ around  
 $z_T(s')/T$ as follows:
\begin{align}
\label{expansion1}
  H^{(0)}\left({\bf s}_{1,m}, \frac {z_{T}(s_m)}T\right)=  
  \sum_{k=0}^{N-1}S^{(k)}\left({\bf s}_{1,m},s',\frac{z_T(s') }{T}\right)+R^{(N)}({\bf s}_{1,m},s'),
\end{align}
where
\begin{align}
&S^{(0)}({\bf s}_{1,m},s',z):=H^{(0)}\left({\bf s}_{1,m},z\right),\label{e:4.36a}\\
&S^{(k)}({\bf s}_{1,m},s',z):=\int_{\Delta_k(s',s_m)}H^{(k)}\left({\bf s}_{1,m},{\bf s}_{1,k}'',z\right)d{\bf s}_{1,k}'',\qquad k=1,2,\ldots \label{SR}\\
&
R^{(N)}({\bf s}_{1,m},s_1'):=\int_{\Delta_N(s',s_m)}H^{(N)}\left({\bf s}_{1,m},{\bf s}_{1,N}'',\frac{z_T(s_N'') }{T}\right)d{\bf s}_{1,N}''.\nonumber
\end{align}

Let us denote
 by $\bbE_{s}$ the conditional expectation with respect to the $\si$-algebra
\begin{equation}
\label{cV}
{\cal V}_s:=\si\left\{V_T(t,z_T(t)),\,t\le s\right\}.
\end{equation}
The following lemma will be used to estimate the remainder term.
{\begin{lemma}
\label{lm031206-18}
Suppose that $m\ge1$, $t_1\ge t_2\ldots \ge t_m\ge0$, $H^{(N)}$, $N\ge0$ are as in \eqref{H0} and  \eqref{e:W_n1x}. Then, 
there exist $\gamma>0$ and $\rho\in(0,1)$ such that for any $N\ge 1$ we can find a
constant $C_N>0$, for which (cf \eqref{010810})
\begin{multline}
\label{021306-18}
\sup_{0\le s_N'\le t_m}\bbE\left\{\E_{s_N'}\left[\int_{\Delta_N({\bf t}_{1,m})}
H^{(N)}\left({\bf s}_{1,m},{\bf s}_{1,N}',\frac{z_T(s_N') }{T}\right)d{\bf s}_{1,m}d{\bf s}_{1,N-1}'\right]\right\}^2\\
\le C_N\left[t_1^{\rho}C^{1-\rho}(t_1) \right]^{2(m+N-1)} T^{2(1-\alpha)-\gamma N}
\end{multline}
for all $  t_1\ge \ldots\ge t_m\ge0$, $T>0$. Here $C(\cdot)$ is defined by
\eqref{021707}.
\end{lemma}}

\bigskip

%
Let ${\frak e}:=(1,\ldots,1)$. It follows  from Proposition \ref{prop011007} that the conditional expectation on the {left} hand side of \eqref{021306-18} has the same distribution as the conditional expectation
\begin{equation*}
 \E_0\left[\int_{\Delta_N({\bf t}_{1,m}-s_N'{\frak e})}
H^{(N)}\left({\bf s}_{1,m},{\bf s}_{1,N-1}',0,0\right)d{\bf s}_{1,m}d{\bf s}_{1,N-1}'\right].
\end{equation*}
Moreover, since ${\cal V}_0\subset \F_0$, where $\F_s$ is defined by \eqref{Fs}, by Jensen's inequality for any random variable $\xi$ we have
\begin{equation*}
 \E\left(\E\left(\xi|{\cal V}_0\right)\right)^2\le\E\left(\E\left(\xi|\F_0\right)\right)^2.
\end{equation*}
Hence, to prove Lemma \ref{lm031206-18}
it suffices to show the following.
\begin{lemma}
\label{lm031206-18a}
Under the assumptions of Lemma $\ref{lm031206-18}$
there exist $\gamma>0$ and $\rho\in(0,1)$ such that for any $N\ge 1$ we can find a constant $C_N>0$, for which
\begin{eqnarray}
\label{021306-18xx}
&&
\bbE\left\{\int_{\Delta_{N-1}({\bf t}_{1,m})}\bbE\left[H^{(N)}\left({\bf s}_{1,m},{\bf s}_{m+1,m+N-1},0,0 \right)\Big|{\cal F}_0\right]d{\bf s}_{1,m}d{\bf s}_{1,N-1}'\right\}^2
\\
&&
\le C_N\left[t_1^{\rho}C^{1-\rho}(t_1) \right]^{2(m+N-1)}T^{2(1-\alpha)-\gamma N}\nonumber
\end{eqnarray}
for all $  t_1\ge \ldots\ge t_{m}\ge 0$ and $T>0$. Here $C(\cdot)$ is
defined by \eqref{021707}.
\end{lemma} 
The proof of this lemma is given in Section \ref{sec5.5-18}.

\bigskip
%
%

Our next result deals with the case when at least one of the factors appearing in the product  considered in \eqref{H0} contains a derivative, that is, belongs to ${\cal D}^{p}(V_T)$
for some $p>0$. 
\begin{lemma}
\label{lm021206-18}
Suppose that $m\ge1$ and $p_j\in\bbZ_+$, 
$H_j\in {\cal D}^{p_j}(V_T)$  (see \eqref{e:DpVT}) and at least
one of the $p_j$-s is positive. Then,
\begin{align}
\label{031206-18}
\lim_{T\to+\infty}\int_{\Delta({\bf t}_{1,m})}\bbE\left[\prod_{j=1}^m 
  H_{j}\left(s_j,0\right)\right]d{\bf s}_{1,m}=0.
\end{align}
\end{lemma}
The proof of the lemma is presented in Section \ref{sec5.4-18}.

\bigskip

To prove \eqref{010306az} we will need a more general result,  which is a consequence of  
Lemmas  \ref{lm031206-18} and \ref{lm021206-18}. For fixed  $m\ge 1$, $n \ge 0$, $p_\ell\in\bbZ_+$, 
$H_\ell\in {\cal D}^{p_\ell}(V_T)$, $\ell\in\bbZ_{m+n}$ and
$t_1\ge\ldots\ge t_{m+n}\ge 0$ let us denote
\begin{equation}
 {I_T^{m,n}( {\bf t}_{1,m+n};(H_\ell))}:=\int_{\Delta( {\bf t}_{1,m+n})}\bbE\left[\prod_{\ell=1}^m
  H_{\ell}\left(s_\ell, \frac{z_T(s_m)}{T}\right)\prod_{\ell=m+1}^{m+n}
  H_{\ell}\left(s_\ell, \frac {z_T(s_\ell)}T\right)\right]d{\bf s}_{1,m+n} \label{e:4.47},
\end{equation}
with the convention that product over an empty set is equal to $1$.
\begin{lemma}
\label{lm011206-18}
{Under the assumptions made in the foregoing we have}
\begin{align}
\lim_{T\to+\infty}\left\{
 I_T^{m,n}( {\bf t}_{1,m+n};(H_\ell))-\int_{\Delta(  {\bf t}_{1,m+n})}\bbE\left[\prod_{\ell=1}^{m+n}
  H_{\ell}\left(s_\ell,0\right)\right]d{\bf s}_{1,m+n}\right\}=0
  \label{010306a}
\end{align}
for any $t_1\ge\ldots\ge t_{m+n}\ge 0$.
\end{lemma}
\proof {The proof is by induction on $n$. If $n=0$, then
\eqref{010306a} follows for an arbitrary $m\in\bbZ_+$ by Proposition
\ref{prop011007} and stationarity of $H_1, \ldots, H_m$ in the time
variable, which in turn follows from stationarity of $V_T$.}

Now suppose that \eqref{010306a} holds for some $n\ge0$ and arbitrary $m\in\bbZ_+$.
{Suppose that  $t_1\ge\ldots\ge t_{m+n+1}\ge 0$ and 
$H_\ell\in {\cal D}^{p_\ell}(V_T)$, for some $p_\ell\in\bbZ_+$, where $\ell\in\bbZ_{m+n+1}$.} Thanks to \eqref{expansion1} we can write
\begin{align}
I_{T}^{m,n+1}( {\bf t}_{1,m+n+1};(H_\ell))=\sum_{k=0}^{N-1}{\cal S}_k(T)+{\cal R}_N(T),
\end{align}
with $N\ge2$ to be chosen later on and
\begin{align*}
&
{\cal S}_k(T):=\int_{\Delta({\bf t}_{1,m+n+1})}\bbE\left[S^{(k)}\left({\bf s}_{1,m},s_{m+1},\frac{z_T(s_{m+1}) }{T}\right)\prod_{\ell=m+1}^{m+n+1}  H_\ell\left(s_\ell, \frac {z_T(s_\ell)}T\right)\right]d{\bf s}_{1,m+n+1}\nonumber\\
&
{\cal R}_N(T):=\int_{\Delta({\bf t}_{1,m+n+1})}\bbE\left[\bbE_{s_{m+1}}\left[R^{(N)}\left({\bf s}_{1,m},s_{m+1}\right)\right]\prod_{\ell=m+1}^{m+n+1}  H_\ell\left(s_\ell, \frac {z_T(s_\ell)}T\right)\right]d{\bf s}_{1,m+n+1}.
\end{align*}
{By \eqref{SR} and \eqref{e:W_n1x} expression ${\cal S}_k(T)$ can be
written as  a sum of terms of the form $I^{m_k,n}({\bf
  t}_{m_k+n};(\tilde H_\ell))$
corresponding to some $m_k\ge 1$, 
$\tilde H_\ell\in {\cal D}^{\tilde p_\ell}(V_T)$ for some $\tilde
p_\ell\in\bbZ_+$ and $\ell\in\bbZ_{
m_k+n}$. When $k\ge 1$ we have to have at least one
 $\tilde p_\ell>0$.
Using the induction hypothesis  and 
Lemma \ref{lm021206-18} we conclude that
\begin{equation*}
\lim_{T\to \infty} {\cal S}_k(T) =0, \qquad \textrm{for}\ k= 1,\ldots,N-1.
\end{equation*}
The only term that may possibly have non-zero limit is 
$$
{\cal S}_0(T)=\int_{\Delta({\bf
    t}_{1,m+n+1})}\bbE\left[\prod_{\ell=1}^{m+1}   H_\ell\left(s_\ell, \frac {z_T(s_{m+1})}T\right)\prod_{\ell=m+2}^{m+n+1}  H_\ell\left(s_\ell, \frac {z_T(s_\ell)}T\right)\right]d{\bf s}_{1,m+n+1}.
$$ 
We can again use the induction hypothesis
 and conclude that
$$
\lim_{T\to+\infty}\left\{{\cal S}_0(T) -\int_{\Delta(  {\bf t}_{1,m+n+1})}\bbE\left[\prod_{\ell=1}^{m+n+1}
  H_{\ell}\left(s_\ell,0\right)\right]d{\bf s}_{1,m+n+1}\right\}=0.
$$ 
Finally, we claim that
\begin{equation}
 \label{011306-18}
 \lim_{T\to \infty}{\cal R}_N(T)=0,
\end{equation}
provided that $N$ is sufficiently large.
Indeed, note that by Proposition \ref{prop011007} the laws of  $
H_\ell\left(s_\ell,z_T(s_\ell)/T\right)$ and that of
$H_\ell\left(0,0\right)$ coincide.  Therefore, using the H\"older inequality, we
can write that (cf \eqref{010810})
\begin{align*}
&
|{\cal R}_N(T)|\le \int_{\Delta({\bf t}_{m+1,m+n+1})}
\left\{\bbE\left[
\int_{\Delta_{m}({\bf t}_{1,m}, s_{m+1})}\bbE_{s_{m+1}}R^{(N)}
\left({\bf s}_{1,m},s_{m+1}\right)
 d{\bf s}_{1,m}
 \right]^2\right\}^{1/2}d{\bf s}_{m+1,m+n+1}\\
 &\hskip 8cm\times\prod_{\ell=m+1}^{m+n+1}\left[\bbE H_\ell^{2(n+1)}\left(0,0\right)\right]^{1/(2n+2)}.
\end{align*}
By \eqref{e:10.24aaa} 
and \eqref{021306-18}
we conclude that there exist $\gamma>0$ and $\rho\in(0,1)$ such that for each $N\ge1$ we can find $C_N>0$, for which
\begin{align}
\label{RNT}
|{\cal R}_N(T)|\le C_N\left[ t_1^{\rho}C^{1-\rho}(t_1)\right]^{N+m-1}T^{-\gamma N}T^{(1-\alpha)(n+1)},\quad t_1\ge \ldots \ge t_m\ge 0,\,T>0.
\end{align}
Choosing $N$ so large that $N\gamma>(1-\alpha)(n+1)$ we conclude \eqref{011306-18}. 
}

{This finishes the proof of the conclusion of  Lemma \ref{lm011206-18} for $n+1$ and any $m\in \Z_+$.}

\qed

The identity \eqref{010306az} is a particular case of Lemma  \ref{lm011206-18}, hence it follows immediately, thus finishing the proof of convergence of moments.\qed

\subsection{Proof of Lemma \ref{lm031206-18a}}

\label{sec5.5-18}

 By \eqref{e:W_n1x} we have 
$$H^{(N)}({\bf s}_{1,m}, {\bf s}_{m+1,m+N-1},0,0)=\sum_{\ell=1}^d\frac 1T\partial_{x_\ell}H^{(N-1)}({\bf s}_{1,m}, {\bf s}_{m+1,m+N-1},0)V_{T,\ell}(0,0).$$
Using the fact that  $V_{T,\ell}(0,0)$ is $\F_0$ measurable and
applying the Cauchy-Schwarz inequality we have
\begin{multline*}
 \bbE\left\{\int_{\Delta_{N-1}({\bf t}_{1,m})}\bbE\left[H^{(N)}\left({\bf s}_{1,m+N-1},0,0 \right)\Big|{\cal F}_0\right]d{\bf s}_{1,m+N-1}\right\}^2\\
\le {d}\sum_{\ell=1}^d\left\{\bbE\left\{\int_{\Delta_{N-1}({\bf t}_{1,m})}\bbE\left[\frac 1T \partial_{x_\ell}H^{(N-1)}\left({\bf s}_{1,m+N-1},0 \right)\Big|{\cal F}_0\right]d{\bf s}_{1,m}d{\bf s}_{1,N-1}'\right\}^4\right\}^{1/2}
\left\{\bbE V_{T,\ell}^4(0,0)\right\}^{1/2}
\end{multline*}
Since  both the random variables $V_{T,\ell}(0,0)$ and conditional
expectation appearing above have a finite Wiener chaos expansion
(see Proposition \ref{lem:A3}) their  $L^4$
norms are equivalent with the $L^2$ norms, see Theorem 3.50, p. 39 of \cite{janson}.
In view of \eqref{e:10.24aa}
and writing $N$ instead of $N-1$,
to prove the lemma it suffices to show that there exist $\ga>0$,
$\rho\in(0,1)$ 
such that
\begin{equation}
 \bbE\left\{\int_{\Delta_{N}({\bf t}_{1,m})}\bbE\left[\frac 1T \partial_{x_\ell}H^{(N)}\left({\bf s}_{1,m+N},0 \right)\Big|{\cal F}_0\right]d{\bf s}_{1,m+N}\right\}^2
 \preceq [t_1^{\rho}C^{1-\rho}(t_1)]^{2(m+N)} T^{-\gamma (N+1)},
\end{equation}
for $t_1\ge \ldots \ge t_m\ge0$ and $T>0$.

To simplify the
notation we  let $H^{(0)}({\bf s}_{1,m},z)$, defined in \eqref{H0}, be
of the form
\begin{equation}
\label{H01}
H^{(0)}\left({\bf s}_{1,m},z\right):=\prod_{j=1}^m
  V_{T,\ell_j}\left(s_j,z\right),
\end{equation}
with $\ell_j\in\bbZ_d$ for $j=1,\ldots,m$. 
To obtain the estimate of the lemma we will use the fact that
$a(\abs{k})$ appearing in the kernel defining $V_T(\cdot)$ is
bounded. In the  general case of  $H_j\in {\cal D}^{p_j}(V_T)$,
recalling \eqref{e:dVj}, we can instead apply the fact that
$\abs{k}^{n}a(\abs{k})$ is bounded for any $n$, since $a$ is
continuous and has compact support. Otherwise the respective estimates are obtained in the same way.

{Clearly, $T^{-1} \partial_{x_\ell}H^{(N)}$ 
is again a product of terms from $\bigcup_{p\in\bbZ_+} {\cal D}^p(V_T)$ .
We can decompose each $V_{T,j}(\cdot)$ appearing in $T^{-1} \partial
{x_\ell}H^{(N)}$ into a sum of $V_{t,j,j'}(\cdot)$, according to
\eqref{VT}, and
represent $H^{(N)}$ as a sum of corresponding products and then estimate each of the terms in this sum separately.}

More precisely, for any multi-index ${\bf u}=(u_1,\ldots, u_m)\in\bbZ_d^m$, let
\begin{equation*}
 H^{(0)}_{{\bf u}}({\bf s}_{1,m},z):=\prod_{i=1}^mV_{T,\ell_i,u_i}(s_i,z).
\end{equation*}
Obviously,
\begin{equation*}
 H^{(0)}({\bf s}_{1,m},z)=\sum_{{\bf u}}H^{(0)}_{{\bf u}}({\bf s}_{1,m},z):=\prod_{i=1}^mV_{T,\ell_i,u_i}(s_i,z),
\end{equation*}
where the summation extends over all indices ${\bf u}$. Moreover, for a multi-index
 ${\bf
  j}:=(j,j',j'')\in\bbZ_d^3$ and ${\bf u}$ as above we let
\begin{equation}
 \W_0^{{\bf u,j}}({\bf s}_{1,m},z):=\frac 1{T}\partial_{
   z_{j''}}H_{{\bf u}}^{(0)}({\bf s}_{1,m},z).
\label{U0}
   \end{equation}
{The indices $j,j'$ are redundant but we have added
them to maintain consistency with the ensuing notation.} 
%
Given $N+2$ multi-indices 
 ${\bf
  j}_i:=(j_i,j'_i,j_i'')\in \bbZ_d^3$, $i=0,\ldots,N+1$ and a multi-index $\bf u$, as above we let 
\begin{equation}
\label{U}
\W_{N+1,T}^{{\bf u,j}_0,\ldots,{\bf j}_{N+1}}({\bf s}_{1,m+N+1},z):=\frac 1{T}\partial_{
   z_{j''_{N+1}}}\left\{\vphantom{\int_0^1}\W_{N,T}^{{\bf
       u,j}_0,\ldots,{\bf j}_{N}}({\bf s}_{1,m+N},z)V_{T,j_{N+1},j_{N+1}'}(s_{m+N+1},z)\right\}.
\end{equation}
{ By virtue of \eqref{U0} we have
\begin{align*}
&\sum_{{\bf u},j_0'', j_1,j_{1}'}\W_{1,T}^{{\bf u},{\bf  j}_0,{\bf
    j}_{1}}({\bf s}_{1,m+1},z)\delta_{j_1,j_{0}''}
=\frac 1T \partial_{z_{j_1''}}\left\{\sum_{{\bf u},j_0'', ,j_{1}'}\vphantom{\int_0^1}\W_{0,T}^{{\bf
       u,j}_0}({\bf s}_{1,m},z)V_{T,j_{0}'',j_{1}'}(s_{m+1},z)\right\}\\
&
=\frac 1T \partial_{z_{j_1''}}H^{(1)}({\bf s}_{1,m+1},z).
\end{align*}}
By induction we can extend the above formula to an arbitrary $N$ and
 obtain
\begin{equation}
\label{UU}
\sum_{{\bf u},j_0'',{\bf j}_1,\ldots,{\bf j}_{N-1},, j_N,j_{N}'}\W_{N,T}^{{\bf u, j}_0,\ldots,{\bf
    j}_{N}}({\bf s}_{1,m+N},z)\prod_{i=1}^{N}\delta_{j_i,j_{i-1}''}
=\frac 1T {\partial_{z_{j_N''}}}H^{(N)}({\bf s}_{1,m+N},z),\quad N\ge 1.
\end{equation}
By \eqref{UU},
to prove \eqref{021306-18xx} it suffices to show that  (cf
\eqref{010810}) there exist $\ga>0$ and $\rho\in(0,1)$ such that for each $N\ge1$
\begin{equation}
\label{021306-18xy}
\bbE\left\{\int_{\Delta_{N}({\bf t}_{1,m})}\bbE\left[\W_{N,T}^{{\bf
       u, j}_0,\ldots,{\bf j}_{N}}\left({\bf s}_{1,m+N},0\right)\Big|{\cal F}_0\right]d{\bf s}_{1,m+N}\right\}^2
 \preceq   [t_1^{\rho}C^{1-\rho}(t_1)]^{2(m+N-1)} T^{-\gamma (N+1)}
\end{equation}
  for arbitrary multi-indices ${\bf
j}_1,\ldots,{\bf j}_{N}\in\bbZ_d^3$, $t_1\ge \ldots t_{m}\ge 0$, $T>0$. 

To simplify the notation, we consider only the case when
{$j_i'=u_\ell=j_0'$, for all $0\le i\le N$ and $1\le \ell\le m$}. Then, the fields appearing in the definition of $\W_{N,T}^{{\bf
       u,j}_0,\ldots,{\bf j}_{N}}(\cdot)$ are based on the same noise. The case when noises may be
independent  can be treated in a similar way. In fact it leads to
better estimates, due to the fact that terms corresponding to
covariance of independent noises vanish. 
To simplify the notation even further we will suppress writing the
multi-indices ${\bf u}$, ${\bf j}_0,\ldots,{\bf j}_{N}$ and simply write $\W_{N,T}$
instead of  $\W_{N,T}^{{\bf u, j}_0,\ldots,{\bf j}_{N}}$. 

Using \eqref{U} and Proposition \ref{lem:A2} we obtain
\begin{equation}
 \W_{N,T}({\bf s}_{1,m+N},z)=\frac 1{T^{N+1}}\sum_{G\in
   \D_{N+m}^2}I_G\left(f(\cdot;{\bf s}_{1,m+N},z)\right),
 \label{e:11.25}
\end{equation}
where $z\in\bbR^d$ and 
\begin{align}
\label{fa}
&f({\bf s}',{\bf k};{\bf s}_{1,m+N},z)
: =\ii^{N+1}\prod_{n=0}^N\left(\sum_{n'=1}^{m+n}( k_{n',1}+  k_{n',2})\right)_{j_n''} \exp\left\{\ii z\cdot \left(\sum_{n=1}^{m+N}( k_{n,1}+
   k_{n,2})\right)\right\}\\
&\times \prod_{n=1}^{m+N}\left[\vphantom{\int_0^1}
 E_T(s_n- s_{n,1}',s_n-  s_{n,2}',  k_{n,1},  k_{n,2})
\ind_{\{ s_{n,1}',  s_{n,2}'\le
   s_n\}}
 \right]\prod_{n=1}^{m+N}\Gamma_{j_n,j'_0}( k_{n,1}+ k_{n,2}).\notag
 \end{align}
We have used the shorthand notation ${\bf s}'$, ${\bf k}$ for
the ensembles of the respective variables  $s_{n,j}'$ and $k_{n,j}$.
{Symbol $I_G(\cdot)$ denotes the multiple stochastic integral
introduced in \eqref{e:A5}, we omit writing the multi-index in this
case, as all the noises in our situation are identical.}

For a given  $G\in \D_{N+m}^2$
we let
\begin{equation}
\label{JG}
 J_G({\bf t}_{1,m}):=\frac 1{T^{N+1}}\int_{\Delta_N({\bf t}_{1,m})} \E\left[ I_G\left(f(\cdot;{\bf s}_{1,m+N} ,0)\right) \Big|{\cal F}_0\right]d{\bf s}_{1,m+N} .
\end{equation} 
Then, obviously
\begin{equation}
\label{UG}
 \int_{\Delta_N({\bf t}_{1,m})} \E\left[\W_{N,T}({\bf s}_{1,m+N},0)\Big|{\cal F}_0\right]= \sum_{G\in
   \D_{N+m}^2} J_G({\bf t}_{1,m}).
\end{equation}


Observe that if $\left((\ell,j),(\ell',j')\right)\in G_{links}$, with $\ell<\ell'$,
then $s_{\ell}>s_{\ell'}$ and
\begin{align}
\label{011707}
& \left[\frac {r_T{\cal E}_T(| k_{\ell,j}|) r_T{\cal E}_T(|| k_{\ell',j'}|)}{(| k_{\ell,j}|| k_{\ell',j'}|)^{d-1}}\right]^{1/2}\int_{\R^2}\exp\left\{-\frac
   12r_T(| k_{\ell,j}|)(s_{\ell}- s_{\ell,j}')-\frac
   12r_T(| k_{\ell',j'}|)(s_{\ell'}- s_{\ell',j}')\right\} \nonumber\\
&
\\
&
\times { \ind_{[s_{\ell,j}'\le s_{\ell}]} \ind_{[ s_{\ell',j'}'\le
  s_{\ell'}]}}\delta( s_{\ell,j}'- s_{\ell',j'}')\delta( k_{\ell,j}+ k_{\ell',j'})d s_{\ell,j}'d s_{\ell',j'}'
 =\tilde  e_T(s_{\ell}-s_{\ell'}, k_{\ell,j})\delta( k_{\ell,j}+ k_{\ell',rj}),\nonumber
\end{align}
with
\begin{equation}
\label{te}
\tilde e_T(s,k):=\frac {e^{-\frac12r_T(|k|)|s|}a(|k|/T) }{|k|^{d-\alpha-1}},\quad (s,k)\in\bbR^{1+d},\,T>0.
\end{equation}
From \eqref{e:A5}, Proposition \ref{lem:A3}  and \eqref{011707} we obtain
\begin{align}
 \label{e:A5a}
& \frac 1{T^{N/(2\beta)}}\int_{\Delta_N({\bf t}_{1,m})} \E\left[ I_G\left(f(\cdot;{\bf s}_{1,m+N} ,0)\right) \Big|{\cal F}_0\right]d{\bf s}_{1,m+N}
\\
 &
 =\int_{\bbR^{(1+d)f(G)}} \left\{\int_{\Delta_N({\bf t}_{1,m})}  H_G\left({\bf s}_{1,m+N},{\bf s}_{G_{free}}',{\bf k}_{G_{free}}\right)  d{\bf s}_{1,m+N} \right\}
\prod_{(\ell,j)\in G_{free}}W_{j'}(d s_{\ell,j}',d k_{\ell,j}),\nonumber
\end{align}
where ${\bf s}_{G_{free}}':=\{s_{\ell,j}';(\ell,j)\in G_{free}\}$ and ${\bf
  k}_{G_{free}}=\{k_{\ell,j};(\ell,j)\in G_{free}\}$ and 
\begin{align}
\label{H}
&
H_G\left({\bf s}_{1,m+N},{\bf s}_{G_{free}}',{\bf
    k}_{G_{free}}\right):=\frac{\ii^{N+1}}{T^{N+1}}\int_{\bbR^{2(1+d)\ell(G)}}
\prod_{n=0}^N\left(\sum_{\ell=1}^{m+n}(k_{\ell,1}+k_{\ell,2})\right)_{j_n''}
  \notag\\
&
\times 
\prod_{(\ell,j)\in
  G_{free}}\left[\vphantom{\int_0^1}e_T(s_\ell-s_{\ell,j}',k_{\ell,j})\ind_{\{s_{\ell,j}'\le
    0\}}\right]\prod_{\left((\ell,j),(\ell',j')\right)\in G_{links}}
 \tilde e_T(s_\ell-s_{\ell'},k_{\ell,j})
\notag
\\
 &
\times \prod_{\ell=1}^{m+N}\Gamma_{j_\ell,j'_0}(k_{\ell,1}+k_{\ell,2})\prod_{\left((\ell,j),(\ell',j')\right)\in
  G_{links}}\delta(k_{\ell,j}+k_{\ell',j'})dk_{\ell,j}dk_{\ell',j'}
\end{align}
and
$$
e_T(s,k):=\left(\vphantom{\int_0^1}r_T(\abs{k})\tilde
  e_T(s,k)\right)^{1/2}
,\quad (s,k)\in\bbR^{1+d},\,T>0.
$$

In what follows we adopt the convention $s_{m+N+1}:=0$ and use the
notation $\tau_\ell:=s_{\ell}-s_{\ell+1}$, $\ell=1,\ldots,m+N$.


 If $(\ell,j)\in G_{free}$, then we can write
\begin{equation}
\label{031707}
 e^{-\frac 12 s_{\ell}r_T(|k_{\ell,j}|)}=\exp\left\{-\frac 12 \sum_{p=\ell}^{m+N}\tau_{p}r_T(|k_{\ell,j}|)\right\}.
\end{equation}
If,{ on the other hand,}  $\{(\ell,j),(\ell',j')\}\in G_{links}$, then
\begin{equation}
\label{041707}
   e^{-\frac 12 (s_{\ell}-s_{\ell'})r_T(|k_{\ell,j}|)}=\exp\left\{-\frac 12 \sum_{p=\ell}^{\ell'-1}\tau_{p}r_T(|k_{\ell,j}|)\right\}.
\end{equation}
Using
\eqref{031707}, \eqref{041707}
we can write
\begin{align}
\label{011807}
 &\prod_{(\ell,j)\in
  G_{free}}\left[\vphantom{\int_0^1}\exp\left\{-\frac12(s_\ell-s_{\ell,j}') r_T(|k_{\ell,j}|)\right\}\ind_{\{s_{\ell,j}'\le
    0\}}\right]\prod_{\left((\ell,j),(\ell',j')\right)\in
  G_{links}}\exp\left\{-\frac12(s_\ell-s_{\ell'})r_T(|k_{\ell,j}|)\right\}\nonumber\\
&
=\exp\left\{-\frac12\sum_{m+N\ge i\ge\ell\ge
    1}\sum_{j=1}^2\tau_i \si_{\ell,j}^i
  r_T(|k_{\ell,j}|)\right\}
\prod_{(\ell,j)\in
  G_{free}}\left[\vphantom{\int_0^1}\exp\left\{\frac12 s_{\ell,j}' r_T(|k_{\ell,j}|)\right\}\ind_{\{s_{\ell,j}'\le
    0\}}\right],
\end{align}
where 
\begin{equation*}
 \sigma_{\ell,j}^{i}=\begin{cases}
  1\ &\textrm{if}\ (\ell,j)\in G_{free},\\
  1  & \textrm{if}\ \left((\ell,j),(\ell',j')\right)\in G_{links}\  \
  \textrm{and }\ell\le i\le \ell'-1,\\
  0& \textrm{otherwise.}
 \end{cases}
\end{equation*}
Fix $i\in\{1,\ldots,N\}$.
Note that
\begin{align}
\prod_{i=0}^{N}&\abs{\sum_{\ell=1}^{m+i}\sum_{j=1}^2k_{\ell,j}}\prod_{\left((\ell,j),(\ell',j')\right)\in
  G_{links}}\delta(k_{\ell,j}+k_{\ell',j'})\notag\\
&=\prod_{i=0}^{N}
\abs{\sum_{\ell=1}^{m+i} \sum_{j=1}^2 \si^i_{\ell,j}k_{\ell,j}}\prod_{\left((\ell,j),(\ell',j')\right)\in
  G_{links}}\delta(k_{\ell,j}+k_{\ell',j'})\notag\\
&\le  \prod_{i=0}^{N}\left(\sum_{\ell=1}^{m+i}\sum_{j=1}^2\sigma_{\ell,j}^{i}|k_{\ell,j}|\right)
\prod_{\left((\ell,j),(\ell',j')\right)\in
  G_{links}}\delta(k_{\ell,j}+k_{\ell',j'})
\label{e:11.28}
\end{align}
Using \eqref{rT}, \eqref{Erb} and \eqref{022806a} we conclude that for some $c_*>0$ we have
\begin{equation}
\label{tep}
\tilde e_T(s,k)\preceq \frac {e^{-c_*|k|^{2\beta}|s|}a(|k|/T) }{|k|^{d+\al-1}}
\end{equation}
and similarly
\begin{equation}
\label{tep1}
e_T(s,k)\preceq \frac {e^{-c_*|k|^{2\beta}|s|}a^{1/2}(|k|/T)}{|k|^{(d+\al-1-2\beta)/2}},\quad (s,k)\in\bbR^{1+d},\,T>0.
\end{equation}
From the above and  \eqref{H} we conclude that
\begin{align}
 &\left|\int_{\Delta_N({\bf t}_{1,m})} H_G\left({\bf s}_{1,m+N},{\bf
       s}_{G_{free}}',{\bf k}_{G_{free}}\right)d{\bf s}_{1,m+N}\right| 
\notag \\
 &
 \preceq T^{-(N+1)}\int_{\bbR^{2(1+d)\ell(G)}} \left\{\prod_{i=0}^{N}\left(\sum_{\ell=1}^{m+i}\sum_{j=1}^2\sigma_{\ell,j}^{i}|k_{\ell,j}|\right)\right\}\nonumber\\
&
\times\left\{\int_{\tilde \Delta_{m+N}(t_1)}\exp\left\{-c_*\sum_{m+N\ge i\ge\ell\ge
    1}\sum_{j=1}^2\tau_i \si_{\ell,j}^i
  |k_{\ell,j}|^{2\beta}\right\}d\tau_{1,m+N}\right\}
\label{aaa}\\
&
\times \prod_{(\ell,j)\in
  G_{free}}\left[\vphantom{\int_0^1}\exp\left\{\frac12 s_{\ell,j}' r_T(|k_{\ell,j}|)\right\}\frac{\ind_{\{s_{\ell,j}'\le
    0\}}}{|k_{\ell,j}|^{(\al-2\beta+d-1)/2}}\right]\notag\\
  &
\times \prod_{\left((\ell,j),(\ell',j')\right)\in
  G_{links}}  \delta(k_{\ell,j}+k_{\ell',j'})\frac{a(\abs{k_{\ell,j}}/T)dk_{\ell,j}dk_{\ell',j'}}{|k_{\ell,j}|^{\al+d-1}},\nonumber
\end{align}
where
$d\tau_{m,n}:=d\tau_\ell\ldots d\tau_n$ for $0\le m\le n$,
$$
\tilde{\Delta}_n(t):=\left[(\tau_1,\ldots,\tau_n):\,\tau_1,\ldots,\tau_n\ge0,\,\sum_{j=1}^n\tau_j\le t\right].
$$
{We can estimate the  integral over the simplex
 appearing in the right hand
side of   (\ref{aaa})   by
\begin{equation}
\label{020813}
\prod_{i=1}^{m+N}\int_0^{t_1}\exp\left\{-c_*\tau_i\sum_{\ell=1}^i\sum_{j=1}^2 \si_{\ell,j}^i
  |k_{\ell,j}|^{2\beta}\right\}d\tau_{i} .
\end{equation}
Using \eqref{021707} and an elementary estimate $(1-e^{-\ga t})/\ga\le
t$ for $\ga,t>0$,
we conclude that for any $\rho\in(0,1)$ we have
\begin{equation}
\label{010813}
\frac{1-e^{-\ga t}}{\ga}\le \frac{t^{\rho}C^{1-\rho}(t)}{(1+\ga)^{1-\rho}},\quad \ga,t>0, 
\end{equation}
with $C(t)$ given by \eqref{021707}.
Performing  the integration over $\tau_i$ in  \eqref{020813} and using
estimate \eqref{010813} we conclude that
\begin{align*}
 & C[t_1^\rho C^{1-\rho}(t_1)]^{m+N}\prod_{i=1}^{m+N}\left(1+\sum_{\ell=1}^i\sum_{j=1}^2\sigma_{\ell,j}^{i}\abs{k_{\ell,j}}  \right)^{-2\beta(1-\rho)}
\end{align*}
with $C$ some constant independent of $t_1,T$.}

Let
$\gamma\in(0,1)$ be arbitrary. Since  $a(\cdot)$ is compactly supported we have
\begin{equation}
\label{KN}
K_{m,N}:=\sup_{\xi_{\ell,j}>0}\prod_{i=0}^{N}\left(\sum_{\ell=1}^{m+i}\sum_{j=1}^2\sigma_{\ell,j}^{i}\abs{\xi_{\ell,j}}\right)^{1-\gamma}\left\{\prod_{\ell=1}^{m+N}
  [a(\abs{\xi_{\ell,1}})a(\abs{\xi_{\ell,2}})]\right\}^{1/2}<+\infty,\quad m,\,N\ge1.
\end{equation}
Using the above  we get
\begin{eqnarray}
 &&\left|\int_{\Delta_{m+N}(t_1)} H_G\left({\bf s}_{1,m+N},{\bf
       s}_{G_{free}}',{\bf k}_{G_{free}}\right)d{\bf s}_{1,m+N}\right|
\nonumber\\
&&
\preceq 
 K_{m,N}
[t_1^\rho C^{1-\rho}(t_1)]^{m+N} T^{-\gamma(N+1)}\int_{\bbR^{2(1+d)\ell(G)}}\left\{\prod_{i=1}^{m+N}\left(1+\sum_{\ell=1}^{i}\sum_{i=1}^2\sigma_{\ell,j}^{i}|k_{\ell,j}|\right)\right\}^{\gamma-2\beta(1-\rho)}\nonumber\\
&&
\label{zzz1}\\
&&
\times \prod_{(\ell,j)\in
  G_{free}}\left[\vphantom{\int_0^1}\exp\left\{c_* s_{\ell,j}'|k_{\ell,j}|^{2\beta}\right\}\frac{\ind_{\{s_{\ell,j}'\le
    0\}}}{|k_{\ell,j}|^{(\al-2\beta+d-1)/2}}\right]\prod_{\left((\ell,j),(\ell',j')\right)\in
  G_{links}}\delta(k_{\ell,j}+k_{\ell',j'})\frac{dk_{\ell,j}dk_{\ell',j'}}{|k_{\ell,j}|^{\al+d-1}},\nonumber
\end{eqnarray}
Now suppose that $\gamma<2\beta(1-\rho)$. Using the fact that  $
\sigma_{\ell,j}^{\ell}$ equals $1$ for all $(\ell,j)\in G_{free}$ and all left
vertices $(\ell,j)$ we obtain
\begin{align*}
&\left\{\prod_{i=1}^{m+N}\left(1+\sum_{\ell=1}^{i}\sum_{j=1}^2\sigma_{\ell,j}^{i}|k_{\ell,j}|\right)\right\}^{\gamma-2\beta(1-\rho)}\prod_{\left((\ell,j),(\ell',j')\right)\in
  G_{links}}\delta(k_{\ell,j}+k_{\ell',j'})\\
&\le \left\{\prod_{\ell=1}^{m+N}\left(1+\sum_{j=1}^2\sigma_{\ell,j}^\ell|k_{\ell,j}|\right)\right\}^{\gamma-2\beta(1-\rho)}\prod_{\left((\ell,j),(\ell',j')\right)\in
  G_{links}}\delta(k_{\ell,j}+k_{\ell',j})\\
&
\preceq \prod_{(\ell,j)\in
  G_{free}}(1+|k_{\ell,j}|)^{\gamma/2-\beta(1-\rho)}\prod_{\left((\ell,j),(\ell',j')\right)\in
  G_{links}}\left[(1+|k_{\ell,j}|)^{\gamma/2-\beta(1-\rho)}\delta(k_{\ell,j}+k_{\ell',j'})\right].
\end{align*}


Therefore,
\begin{eqnarray}
 &&\left|\int_{\Delta_{N}({\bf t}_{1,m})} H_G\left({\bf s}_{1,m+N},{\bf
       s}_{G_{free}}',{\bf k}_{G_{free}}\right)d{\bf s}_{1,m+N}\right|
\nonumber\\
&&
\preceq
[t_1^\rho C^{1-\rho}(t_1)]^{m+N} T^{-(N+1)\gamma}\left\{
\int_{\bbR^d}\frac{
  dk}{|k|^{\al+d-1}(1+|k|)^{\beta(1-\rho)-\gamma/2}}\right\}^{\ell(G)}\nonumber\\
&&
\times\prod_{(\ell,j)\in
  G_{free}}\left[\vphantom{\int_0^1}\exp\left\{c_* s_{\ell,j}' |k_{\ell,j}|^{2\beta}\right\}\frac{(1+|k_{\ell,j}|)^{\gamma/2-\beta(1-\rho)}\ind_{\{s_{\ell,j}'\le
    0\}}}{|k_{\ell,j}|^{(\al-2\beta+d-1)/2}}\right].
\label{zzz3}
\end{eqnarray}
From the definition of $J_G({\bf t}_{1,m})$, see \eqref{JG}, and \eqref{zzz3}
we conclude that
\begin{align}
&\bbE J_G^2({\bf t}_{1,m})= \int_{\bbR^{f(G)(1+d)}}\left[\int_{\Delta_{N}({\bf t}_{1,m})} H_G\left({\bf s}_{1,m+N},{\bf
       s}_{G_{free}}',{\bf k}_{G_{free}}\right)d{\bf s}_{1,m+N}\right]^2 \prod_{(\ell,j)\in
  G_{free}} ds_{\ell,j}'dk_{\ell,j}\nonumber
\\
&
\preceq
 [t_1^\rho C^{1-\rho}(t_1)]^{2(m+N)}T^{-2(N+1)\gamma }\left\{
\int_{\bbR^d}\frac{
  dk}{|k|^{\al+d-1}(1+|k|)^{\beta(1-\rho)-\gamma/2}}\right\}^{2\ell(G)}\\
&
\hspace{3in}\times\left\{\int_{\bbR^d} \frac{dk}{|k|^{\al+d-1}(1+|k|)^{2\beta(1-\rho)-\gamma}}\right\}^{f(G)}.\nonumber
\label{e:11.37}
\end{align}
Since $\al+2\beta\ge \al+\beta>1$, we can choose  sufficiently small
$\gamma>0$, $\rho\in(0,1)$ so the integrals appearing on the utmost right hand side converge.  In light of \eqref{UU} and \eqref{UG}
this
concludes the proof of  
Lemma \ref{lm031206-18a}.
\qed\qed

\medskip

\subsection{Proof of Lemma \ref{lm021206-18}}

\label{sec5.4-18}


%
%

To simplify the notation we shall prove the lemma for 
$$
H_1(s,z)=\frac{1}{T}(\partial_{z_{j''}}V_{j_1'})\left(s,\frac{z}{T}\right),\quad H_\ell(s,z)=V_{j_\ell'}\left(s,\frac{z}{T}\right),\quad \ell=2,\ldots,m
$$
for some $j'',j_1',\ldots,j_m'\in \bbZ_d$. The general case can be dealt with in the same fashion.

Given $\ell\in\bbZ_d$ we can write
 $H_\ell=\sum_{j=1}^dH_{\ell}^{(j)}$, where $H^{(j)}_\ell$ corresponds to   $Y_j$ component of the noise in \eqref{e:10.1}.
For any multi-index ${\bf j}:=(j_1\ldots, j_{m})\in\bbZ_d^m$ let us denote
\begin{equation*}
L^{\bf j}( {\bf s}_{1,m}):= \bbE\left[\prod_{\ell=1}^m 
  H_{\ell}^{(j_\ell)}\left(s_j,0\right)\right].
\end{equation*}
Then 
\begin{equation}
\label{012707}
 \bbE\left[\prod_{\ell=1}^m 
  H_{\ell}\left(s_j,0\right)\right]=\sum_{{\bf j}}L^{\bf j}( {\bf s}_{1,m})
  \end{equation}
  The proof of the lemma shall be completed as soon as we show that
  \begin{align}
\label{031206-18j}
\lim_{T\to+\infty}\int_{\Delta({\bf t}_{1,m})}\bbE\left[L^{\bf j}( {\bf s}_{1,m})\right]d{\bf s}_{1,m}=0
\end{align}
for each ${\bf j}\in\bbZ_d^m$.
We consider only the case when  {$j_0=\ldots=j_{m}$}, as the  cases corresponding to the other indices can be dealt with similarly. As before we drop writing the multi-index. According to \eqref{e:A7} we have
\begin{equation}
\label{022107}
 {\bbE}\left[L( {\bf s}_{1,m})\right]=\sum_{G\in \G_{m}^2}I_G\left(\tilde f(\cdot; {\bf s}_{1,m})\right)
\end{equation}
where $I_G$ is given by \eqref{e:A5} and 
\begin{align}
\label{f}
&\tilde f({\bf s}',{\bf k},{\bf s}_{1,m})
 =\ii^{m}T^{-1} (k_{1,1}+k_{1,2})_{j''}\notag\\
 &\times \prod_{\ell=1}^{m}\left[\vphantom{\int_0^1}
 E_T(s_\ell-s_{\ell,1}',s_\ell-s_{\ell,2}',k_{\ell,1},k_{\ell,2})
\ind_{\{s_{\ell,1}', s_{\ell,2}'\le s_\ell\}}\Gamma_{j_\ell',j_0}(k_{\ell,1}+k_{\ell,2}).
 \right].
\end{align}
Here the summation $\sum_{G\in \G_{m}^2}$ extends over all complete Feynmann diagrams with $m$ nodes and $2$ hands, $E_T$ is defined in \eqref{ET} and   ${\bf s}'$ and ${\bf k}$ denote the ensembles of variables $\left(s_{\ell,j}'\right)$ and $\left(k_{\ell,j}\right)$, respectively.
We adopt the convention that $s_{m+1}=0$.
Invoking the definition of
$E_T(s_\ell-s_{\ell,1}',s_\ell-s_{\ell,2}',k_{\ell,1},k_{\ell,2})$ and
assumptions \eqref{Erb} -- \eqref{022806a} 
it is clear  that $\tilde f$ converges pointwise to $0$ as $T\to \infty$. Moreover, since $a$ has a bounded support we have
\begin{equation}
\label{010811}
 \abs{\tilde f({\bf s}',{\bf k},{\bf s}_{1,m})}\preceq g({\bf s}',{\bf k},{\bf s}_{1,m}),
\end{equation}
where 
\begin{equation}
g({\bf s}',{\bf k},{\bf s}_{1,m})
=\prod_{\ell=1}^{m}\prod_{j=1}^2
\frac{e^{-c_*\abs{k_{\ell,j}}^{2\beta} (s_\ell-s_{\ell,j}')}\ind_{[s_{\ell,j}'\le s_\ell]}}{\abs{k_{\ell,j}}^{(\alpha+d-1)/2-\beta}},
\label{e:funct_g}
\end{equation}
for some positive constant $c^*$.
Hence, by \eqref{010811}, \eqref{022107} and the dominated convergence theorem to prove \eqref{031206-18j}
it is enough to show that 
\begin{equation}
 \label{e:JG}
 J_G:=\int_0^{t}\int_{\Delta_{m}({\bf t}_{1,m})}I_G(g(\cdot, {\bf s}_{1,m}))d{\bf s}_{1,m}<\infty\quad\mbox{ for any $G\in \G_{m}^2$}. 
\end{equation}


Recalling  the definition \eqref{e:A5}  and integrating out with respect to variables $s_{\ell,j}'$ we obtain that for each $m$ there exists a constant $C>0$, independent of $t_1$, such that
\begin{equation*}
 J_G\le 
 C\int_{\Delta_{m}(t_1)}d{\bf s}_{1,m}\int_{\R^{2md}}\prod_{\substack{((\ell,j),(\ell',j'))\in G\\ \ell<\ell'}} 
 \frac
 {e^{-c_*\abs{k_{\ell,m}}^{2\beta}(s_\ell-s_{\ell'})}}{\abs{k_{\ell,j}}^{\alpha+d-1}}\delta(k_{\ell,j}+k_{\ell',j'})d{\bf
   k},\quad T,t_1>0 .
\end{equation*}
Substituting $\tau_\ell:=s_\ell-s_{\ell+1}$, estimating
$e^{-c_*\abs{k_{\ell,j}}^{2\beta}(s_\ell-s_{\ell'})}\le
e^{-c_*\abs{k_{\ell,j}}^{2\beta}\tau_\ell}$ and enlarging the domain
of integration from the simplex $[\sum_{\ell=1}^m \tau_\ell\le t_1,\,\tau_\ell\ge0,\,\ell=1,\ldots,m]$ to $[0,t_1]^{m}$ we get
\begin{equation}
 J_G\le C \int_{[0,t_1]^{m}}\int_{\R^{2md}}\prod_{\substack{((\ell,j),(\ell',j'))\in G\\ \ell<\ell'}} 
 \frac {e^{-c_*\abs{k_{\ell,j}}^{2\beta}\tau_\ell}}{\abs{k_{\ell,j}}^{\alpha+d-1}}\delta(k_{\ell,j}+k_{\ell',j'})dk_{\ell,j}dk_{\ell',j'} d{\bf \tau}_{1,m}.
\label{e:JG_estimate}
 \end{equation}
When we perform the integration  with respect to $\tau_\ell$ in the
right hand side  there are three possible outcomes:

\noindent {{\em Case 1.} both $(\ell,1)$ and $(\ell,2)$ are left
vertices of some links in $G$. Denote the set of those   $\ell$-s by
$L_1(G)$ and its cardinality by $l_1(G)$. Then, (cf \eqref{021707})
\begin{eqnarray}
&& \int_0^{t_1}
  e^{-c_*(\abs{k_{\ell,1}}^{2\beta}+\abs{k_{\ell,2}}^{2\beta})\tau_\ell}d\tau_\ell=
  \frac{1-e^{-c_*(\abs{k_{\ell,1}}^{2\beta}+\abs{k_{\ell,2}}^{2\beta})t_1}}{c_*(\abs{k_{\ell,1}}^{2\beta}+\abs{k_{\ell,2}}^{2\beta})}.
\label{e:JGa}
 \end{eqnarray}
{\em Case 2.}  only one of $(\ell,1)$ and $(\ell, 2)$ is a left vertex of a link in $G$, say $(\ell,j)$. Denote the set of those   $\ell$-s by
$L_2(G)$ and its cardinality by $l_2(G)$.  Then
\begin{equation}
 \int_0^{t_1} e^{-c_*\abs{k_{\ell ,j}}^{2\beta}\tau_\ell}d\tau_\ell=\frac{1-e^{-c_*\abs{k_{\ell,j}}^{2\beta}t_1}}{c_*\abs{k_{\ell ,j}}^{2\beta}};
\label{e:JGb}
 \end{equation}
{\em Case 3.} none of $(\ell,1)$, $(\ell, 2)$ is a left vertex in
$G$. The cardinality of the set of those $\ell$ obviously equals $l_1(G)$. In this case the integral with respect to $\tau_\ell$ equals
$t_1$. }

\bigskip

{Note an obvious identity 
\begin{equation}
\label{lr}
2l_1(G)+l_2(G)=m.
\end{equation}
Performing the change of variables $t_1^{1/(2\beta)} k_{\ell,j}\mapsto
k_{\ell,j}$, when $(\ell,j)\in {\cal L}(G)$- the set of all left vertices of some links, and using \eqref{lr}  the right
hand side of \eqref{e:JG_estimate}  can we rewritten as
\begin{align}
&\label{010812}
\frac{C t_1^{mH}}{c_*^{l_1(G)+l_2(G)}} \int_{\R^{md}}\prod_{\ell\in L_1(G)}
  \frac{1-e^{-c_*(\abs{k_{\ell,1}}^{2\beta}+\abs{k_{\ell,2}}^{2\beta})}}{\abs{k_{\ell,1}}^{2\beta}+\abs{k_{\ell,2}}^{2\beta}}
\mathop{\prod_{\ell\in L_2(G)}}\limits_{(\ell,j)\in {\cal L}(G)}
\frac{1-e^{-c_*\abs{k_{\ell,j}}^{2\beta}}}{\abs{k_{\ell
  ,j}}^{2\beta}}\prod_{(\ell,j)\in{\cal L}(G)}\frac{
  dk_{\ell,j}}{|k_{\ell,j}|^{\al+d-1}}\nonumber\\
&
\preceq  t_1^{mH} \prod_{\ell\in L_1(G)}\int_{\R^{2d}}
  \frac{1-e^{-c_*(\abs{k_{\ell,1}}^{2\beta}+\abs{k_{\ell,2}}^{2\beta})}}{\abs{k_{\ell,1}}^{2\beta}+\abs{k_{\ell,2}}^{2\beta}}\cdot\frac{
  dk_{\ell,1}dk_{\ell,2}}{|k_{\ell,1}|^{\al+d-1}|k_{\ell,2}|^{\al+d-1}}\\
&
\times\mathop{\prod_{\ell\in L_2(G)}}\limits_{(\ell,j)\in {\cal L}(G)}\int_{\R^{d}}
\frac{1-e^{-c_*\abs{k_{\ell,j}}^{2\beta}}}{\abs{k_{\ell
  ,j}}^{2\beta}}\cdot\frac{
  dk_{\ell,j}}{|k_{\ell,j}|^{\al+d-1}}.\nonumber
\end{align}
Invoking an elementary fact that $(1-e^{-x})/x\approx (1+x)^{-1}$,
$x>0$ and \eqref{e:11.11a} we conclude that all integrals in the right hand side of
\eqref{010812} converge, as $\al+\beta>1$ and $\al<1$. 
The proof of the lemma is therefore finished.\qed}

\bigskip

{\begin{remark}\label{rem:tightness}
 Note that, in light of \eqref{lr}, the above argument shows that 
\begin{equation}
\label{020812}
J_G\preceq t_1^{mH},\quad T,t_1>0.
\end{equation}
This estimate will be useful in the proof of tightness.
\end{remark}}

\subsection{Tightness - the end of the proof of Theorem \ref{thm:main}}

\label{sec:tightness}

In light of the already proved Theorem \ref{thm090801} the only item
that requires yet to be shown is tightness of the laws of
$\left(z_T(t)\right)_{t\ge0}$, $T>0$ over
$C([0,+\infty);\R^d)$.
The processes $\left(z_T(t)\right)_{t\ge0}$, $T>0$ have stationary increments, therefore, to demonstrate  tightness, as $T\to+\infty$, it suffices to show that for any $S>0$ there exist $C>0$ and $\kappa>0$ such that 
\begin{equation}
 \E \abs{z_T(t)}^2\le C t^{1+\kappa}, \qquad t\in[0,S], T>1.
\label{e:12.1}
 \end{equation}
{Thanks to the stationarity of
 $\left(V_T(s,z_T(s)/T)\right)_{s\in\bbR}$ (cf Proposition
 \ref{prop011007})  we have 
 \begin{align*}
 & \E \abs{z_T(t)}^2=2
   \int_0^t\int_u^t\E\left[V_T\left(s,\frac{z_T(s)}T\right)\cdot V_{T}\left(u,\frac {z_T(u)}T\right)\right]duds
  \\
&
=2 \int_0^tdu \int_0^u\E\left[ V_{T}\left(s,\frac{z_T(s)}T\right)\cdot
  V_{T}\left(0,0\right)\right]ds.
 \end{align*}
We expand $V_T(s,z_T(s)/T)$ around $0$, according to
\eqref{expansion1}.  We can easily see by a direct calculation,
similar to those made in Section \ref{sec:4.2}, that 
\begin{equation*}
 \int_0^tdu\int_0^u \E\left[V_T(s,0)\cdot V_T(0,0)\right]ds \preceq t^{2H},\quad
   t>0, \,T>1.
\end{equation*}
By Remark \ref{rem:tightness}, see estimate \eqref{020812}, the terms resulting from $S^{(k)}$,
$k=1,\ldots,N-1$ (cf \eqref{expansion1})  can be bounded from above by $C
t^{(k+1)H}$ for $t\in[0,S]$, where constant $C>0$ is independent of
$T>0$. Considering the remainder term, according to \eqref{RNT}, it can
be estimated by $Ct^{\rho(N+1)}$ for some $\rho\in(0,1)$, independent
of $N$. Choosing the latter so large that $\rho(N+1)>1$ and remebering  that
$(k+1)H\ge2H>1$ for $k=1,\ldots,N-1$ we conclude \eqref{e:12.1},
finishing in this way the proof of tightness and Theorem \ref{thm:main}.\qed}
  
%
%

\section{Proof of Proposition \ref{prop:Rosenblatt}}

\label{sec5}

To simplify we assume that $a_0$, $r_0$ in \eqref{Einfty} both equal $1$.
Since the Rosenblatt process is determined by its moments we can use the moment method to establish the result.

Fix any $N\in \bbN$, $b_1,\ldots, b_N\in \bbR$ and $r_1,\ldots, r_N\in \bbR_+$ and denote 
\begin{equation}
 \psi(t)=\sum_{m=1}^N b_m\ind_{[0,r_m]}(t).
 \label{e:10.16a}
\end{equation}
For a given  $j=1,\ldots,d$ we can write
\begin{align}
 \sum_{m=1}^N b_mZ_j(r_m)=
 &\int_{\bbR^{2d+2}}F(s,k,s',k')
W_j(ds,dk)W_j(ds',dk'),
\label{e:10.17a}
\end{align}
where
\begin{align}
F(s,k,s',k'):=
 &\int_{\bbR}\psi(t)\ind_{(-\infty,t]^2}(s,s')
  \left(\abs{k}\abs{k'}\right)^{-(\al+d-1)/2-\beta}
\exp\left\{-\frac12|k|^{2\beta}(t-s)\right\}
\notag\\
&
\exp\left\{-\frac12|k'|^{2\beta}(t-s')\right\}
dt.
\label{e:10.19}
\end{align}

By Lemma \ref{lem:A2} we have
\begin{equation}
  \E\left(\sum_{m=1}^Nb_mZ_j(r_m)\right)^n=
 \sum_{G\in \G_n^2}I_G,
 \label{e:10.20}
\end{equation}
where
\begin{align}
I_G= 
 \int_{\bbR^{2nd+2n}}&F(s_{1,1},k_{1,1}, s_{1,2},k_{1,2})\ldots
 F(s_{n,1},k_{n,1}, s_{n,2},k_{n,2})\notag\\
 &\prod_{((\ell,m),(\ell',m'))\in G}\delta(s_{\ell,m}-s_{\ell',m'})\delta(k_{\ell,m}+k_{\ell',m'})ds_{\ell,m}ds_{\ell',m'}dk_{\ell,m}dk_{\ell',m'}
 \label{e:10.21}
\end{align}
Clearly if $G=G_1\cup\ldots \cup G_m$, where $G_i$ are connected components of $G$, then $I_G=I_{G_1}\ldots I_{G_m}$ 
Following the proof of Theorem 3.2 in \cite{BGT_Rosenblatt}, it suffices to show that there exists $C>0$ such that for any $G\in \tilde \G_n^2$ (the set of complete diagrams with a single connected component)
\begin{equation}
I_G=C^n\int_{\bbR^n}\prod_{j=1}^n\left[ \psi(t_j)\abs{t_j-t_{j+1}}^{H-1}\right]d{\bf t}_{1,n},\quad n\ge1,
\label{e:10.22}
\end{equation}
with $H$ given by \eqref{e:H} and $t_{n+1}:=t_1$.
If $G\in \tilde  \G_n^2$ then \eqref{e:10.21} can be written as
\begin{equation}
\label{022507}
I_G= 
 \int_{\bbR^{n(d+1)}}F(s_1,k_1, s_2,-k_2)F(s_2, k_2,s_3,-k_3)\ldots
 F(s_n,k_n, s_1,-k_1)d{\bf k}_{1,n}d{\bf s}_{1,n}.
\end{equation}

Note that
\begin{align}
& \int_{s\vee s'}^{+\infty}\exp\left\{-\abs{k}^{2\beta}[t-(s+s')/2]\right\}dt=\frac{e^{-\frac 12\abs{s-s'}\abs{k}^{2\beta}}}{\abs{k}^{2\beta}}.
\label{e:10.6}
 \end{align}
We perform integration over ${\bf s}_{1,n}$ in \eqref{022507}. Using \eqref{e:10.6} we obtain
\begin{align}
 I_G= \int_{\bbR^{n(d+1)}}\prod_{j=1}^n\left[\psi(t_j)
\frac{e^{-\frac 12\abs{k_j}^{2\beta}\abs{t_j-t_{j+1}}}}{\abs{k_j}^{\al+d-1}}
\right]
 d{\bf k}_{1,n}d{\bf t}_{1,n}.
\label{e:10.23}
 \end{align}
Substituting 
\begin{equation}
\label{012607}
k_i':=k_i\abs{t_i-t_{i-1}}^{\frac 1{2\beta}}, \quad i=1,\ldots, n
\end{equation}
 we obtain
\begin{equation*}
  I_G= \left(\int_{\bbR^d} \frac{e^{-\frac 12\abs{k}^{2\beta}}}{\abs{k}^{\al+d-1}}dk\right)^n \int_{\bbR^n}
\prod_{j=1}^n\left[\psi(t_j) \abs{t_i-t_{i+1}}^{H-1}\right]
 d{\bf t}_{1,n},
\end{equation*}
(cf \eqref{e:H}) with the convention $t_{n+1}:=t_1$. Hence \eqref{e:10.22} follows.\qed


 \setcounter{equation}{0}
 \setcounter{thm}{0}
\renewcommand{\thesection}{A}

\end{document}